\documentclass[smallextended]{svjour3}   
\smartqed  
\usepackage{graphicx}
\usepackage{xcolor}

\usepackage{bm,amsfonts,amsmath,amssymb}
\usepackage[ruled,vlined,linesnumbered]{algorithm2e}

\DeclareMathOperator*{\argmin}{argmin}

\usepackage{hyperref} 
\hypersetup{
colorlinks=false, 
breaklinks=true,
urlcolor= blue, 
linkcolor= blue,
citecolor=red, 
}

\newcommand{\sigmab}{\underline{\sigma}}

\journalname{Constraints}
\begin{document}
\title{The Potential of Quantum Annealing for Rapid Solution Structure Identification}
\titlerunning{Quantum Annealing for Rapid Solution Structure Identification}

\author{
Yuchen Pang \and
Carleton Coffrin \and
Andrey Y. Lokhov \and
Marc Vuffray
}


\institute{
Y. Pang, Graduate Research Assistant \at
    University of Illinois at Urbana-Champaign, Department of Computer Science \\
    Champaign, IL, 61801, USA \\
    \email{yuchenp2@illinois.edu} \\
    ORCID: 0000-0002-4532-7053
\and
C. Coffrin, Staff Scientist \at
    Los Alamos National Laboratory, Advanced Network Science Initiative \\
    Los Alamos, NM, 87545, USA \\
    \email{cjc@lanl.gov} \\
    ORCID: 0000-0003-3238-1699
\and
A. Y. Lokhov, Staff Scientist \at
    Los Alamos National Laboratory, Advanced Network Science Initiative \\
    Los Alamos, NM, 87545, USA \\
    \email{lokhov@lanl.gov} \\
    ORCID: 0000-0003-3269-7263
\and
M. Vuffray, Staff Scientist \at
    Los Alamos National Laboratory, Advanced Network Science Initiative \\
    Los Alamos, NM, 87545, USA \\
    \email{vuffray@lanl.gov} \\
    ORCID: 0000-0001-7999-9897
}

\date{}

\maketitle
\begin{abstract}
The recent emergence of novel computational devices, such as quantum computers, coherent Ising machines, and digital annealers presents new opportunities for hardware-accelerated hybrid optimization algorithms.  Unfortunately, demonstrations of unquestionable performance gains leveraging novel hardware platforms have faced significant obstacles.  One key challenge is understanding the algorithmic properties that distinguish such devices from established optimization approaches.  Through the careful design of contrived optimization tasks, this work provides new insights into the computation properties of quantum annealing and suggests that this model has the potential to quickly identify the structure of high-quality solutions.  A meticulous comparison to a variety of algorithms spanning both complete and local search suggests that quantum annealing's performance on the proposed optimization tasks is distinct.  This result provides new insights into the time scales and types of optimization problems where quantum annealing has the potential to provide notable performance gains over established optimization algorithms and suggests the development of hybrid algorithms that combine the best features of quantum annealing and state-of-the-art classical approaches.

\keywords{Discrete Optimization, Ising Model, Quadratic Unconstrained Binary Optimization, Local Search, Quantum Annealing, Large Neighborhood Search, Integer Programming, Belief Propagation}

\end{abstract}

\section{Introduction}

As the challenge of scaling traditional transistor-based Central Processing Unit (CPU) technology continues to increase, experimental physicists and high-tech companies have begun to explore radically different computational technologies, such as quantum computers \cite{ibm_quantum,45919,chmielewski2018cloud}, quantum annealers \cite{PhysRevE.58.5355,johnson2011quantum} and coherent Ising machines \cite{mcmahon2016fully,Inagaki603,7738704}.  The goal of all of these technologies is to leverage the dynamical evolution of a physical system to perform a computation that is challenging to emulate using traditional CPU technology, the most notable example being the simulation of quantum physics \cite{Feynman1982-FEYSPW}.  Despite their entirely disparate physical implementations, optimization of quadratic functions over binary variables (e.g., the Quadratic Unconstrained Binary Optimization (QUBO) and Ising models \cite{RevModPhys.39.883}) has emerged as a challenging computational task that a wide variety of novel hardware platforms can address.  As these technologies mature, it may be possible for this specialized hardware to rapidly solve challenging combinatorial problems, such as Max-Cut \cite{Haribara2016} or Max-Clique \cite{10.3389/fphy.2014.00005}, and preliminary studies have suggested that some classes of Constraint Satisfaction Problems can be effectively encoded in such devices because of their combinatorial structure \cite{10.3389/fict.2016.00014,10.3389/fphy.2014.00056,Rieffel2015,1506.08479}.

At this time, understanding the computational advantage that these hardware platforms may bring to established optimization algorithms remains an open question.  For example, it is unclear if the primary benefit will be dramatically reduced runtimes due to highly specialized hardware implementations \cite{7063111,6662276,fuhitsu_da} or if the behavior of the underlying analog computational model will bring intrinsic algorithmic advantages \cite{albash2018adiabatic,1808.09999}.
A compelling example is gate-based quantum computation (QC), where a significant body of theoretical work has found key computational advantages that exploit quantum properties \cite{shor1994algorithms,grover,coles2018quantum}. 
Indeed, such advantages have recently been demonstrated on quantum computing hardware for the first time \cite{Arute2019}.  Highlighting similar advantages on other computational platforms, both in theory and in practice, remains a central challenge for novel physics-inspired computing models \cite{kalinin2018global,leleu2019destabilization,hamerly2019experimental}.

Focusing on quantum annealing (QA), this work provides new insights on the properties of this computing model and identifies problem structures where it can provide a computational advantage over a broad range of established solution methods.  The central contribution of this work is the analysis of {\em tricky} optimization problems (i.e., Biased Ferromagnets, Frustrated Biased Ferromagnets,  and Corrupted Biased Ferromagnets) that are challenging for established optimization approaches but are easy for QA hardware, such as D-Wave's 2000Q platform.  This result suggests that there are classes of optimization problems where QA can effectively identify global solution structure while established heuristics struggle to escape local minima.  Two auxiliary contributions that resulted from this pursuit are the identification of the Corrupted Biased Ferromagnet problem, which appears to be a useful benchmark problem beyond this particular study, and demonstration of the most significant performance gains of a quantum annealing platform to the established state-of-the-art alternatives, to the best of our knowledge.

This work begins with a brief introduction to both the mathematical foundations of the Ising model, Section \ref{sec:ising}, and quantum annealing, Section \ref{sec:qa_foundation}.  It then reviews a variety of algorithms than can be used to solve such models in Section \ref{sec:algs}.  The primary result of the paper is presented in carefully designed structure detection experiments in Section \ref{sec:structure}.  Open challenges relating to developing hybrid algorithms are discussed in Section \ref{sec:hybrid}, and Section \ref{sec:conclusion} concludes the paper.

\section{A Brief Introduction to Ising Models}
\label{sec:ising}

This section introduces the notations of the paper and provides a brief introduction to Ising models, a core mathematical abstraction of QA.
The Ising model refers to the class of graphical models where the nodes, ${\cal N} = \left\{1,\dots, N\right\}$, represent {\em spin} variables (i.e., $\sigma_i \in \{-1,1\} ~\forall i \in {\cal N}$), and the edges, ${\cal E} \subseteq {\cal N} \times {\cal N}$, represent pairwise {\em interactions} of spin variables (i.e., $\sigma_i \sigma_j ~\forall i,j \in {\cal E}$).  A local {\em field} $\bm h_i ~\forall i \in {\cal N}$ is specified for each node, and an interaction strength $\bm J_{ij} ~\forall i,j \in {\cal E}$ is specified for each edge. The {\em energy} of the Ising model is then defined as:
\begin{align}
    E(\sigma) &= \sum_{i,j \in {\cal E}} \bm J_{ij} \sigma_i \sigma_j + \sum_{i \in {\cal N}} \bm h_i \sigma_i \label{eq:ising_eng}
\end{align}
Originally introduced in statistical physics as a model for describing phase transitions in ferromagnetic materials \cite{gallavotti2013statistical}, the Ising model is currently used in numerous and diverse application fields such as neuroscience \cite{hopfield1982neural,schneidman2006weak}, bio-polymers \cite{morcos2011direct}, gene regulatory networks \cite{marbach2012wisdom}, image segmentation \cite{panjwani1995markov}, statistical learning \cite{NIPS2016_6375,lokhov2018optimal,vuffray2019efficient}, and sociology \cite{eagle2009inferring}.

This work focuses on finding the lowest possible energy of the Ising model, known as a {\em ground state},  that is, finding the globally optimal solution of the following discrete optimization problem:
\begin{align}
    & \min: E(\sigma) \label{eq:ising_opt}\\
    & \mbox{s.t.: } \sigma_i \in \{-1, 1\}  ~\forall i \in {\cal N} \nonumber
\end{align}
The coupling parameters of Ising models are categorized into two groups based on their sign: the \emph{ferromagnetic} interactions $\bm  J_{ij} < 0$, which  encourage neighboring spins to take the same value, i.e., $\sigma_i \sigma_j = 1$, and \emph{anti-ferromagnetic} interactions $\bm  J_{ij}>0$, which encourage neighboring spins to take opposite values, i.e., $\sigma_i \sigma_j = -1$.

\paragraph{Frustration:}
The notion of frustration is central to the study of Ising models and refers to any instance of \eqref{eq:ising_opt} where the optimal solution does not achieve the minimum of all local interactions \cite{ising_frustration}. Namely, the optimal solution of a frustrated Ising model, $\sigma^*$, satisfies the following property:
\begin{align}
    E(\sigma^*) > \sum_{i,j \in {\cal E}} - |\bm J_{ij}| - \sum_{i \in {\cal N}} |\bm h_i|
\end{align}
%
%
%

\paragraph{Gauge Transformations:}
A valuable property of the Ising model is the gauge transformation, which characterizes an equivalence class of Ising models.  Consider the optimal solution of Ising model $S$, $\bm \sigma^{s}$.  One can construct a new Ising model $T$ where the optimal solution is the target state $\bm \sigma^{t}$ by applying the following parameter transformation:
\begin{subequations}
\begin{align}
    \bm J^t_{ij} &= \bm J^s_{ij} \bm \sigma^s_i \bm \sigma^s_j \bm \sigma^t_i \bm \sigma^t_j ~\forall i,j \in {\cal E} \\
    \bm h^t_i &= \bm h^s_i \bm \sigma^s_i \bm \sigma^t_i ~\forall i \in {\cal N}
\end{align}
\end{subequations}
This $S$-to-$T$ manipulation is referred to as a gauge transformation.   Using this property, one can consider the class of Ising models where the optimal solution is $\sigma_i = -1 ~\forall i \in {\cal N}$ or any arbitrary vector of ${-1,1}$ values without loss of generality.

\paragraph{Classes of Ising Models:}
Ising models are often categorized by the properties of their optimal solutions with two notable categories being Ferromagnets (FM) and Spin glasses.
Ferromagnetic Ising models are unfrustrated models possessing one or two optimal solutions.  The traditional FM model is obtained by setting $\bm J_{ij}= -1, \bm h_i = 0$. The optimal solutions have a structure with all spins {\em pointing} in the same direction, i.e., $\sigma_i = 1$ or $\sigma_i = -1$, which mimics the behavior of physical magnets at low temperatures.
%
In contrast to FMs, Spin glasses are highly frustrated systems that exhibit an intricate geometry of optimal solutions that tend to take the form of a hierarchy of isosceles sets \cite{mezard1985microstructure}.
Spin glasses are challenging for greedy and local search algorithms \cite{barahona1982computational} due to the nature of their energy landscape \cite{mezard2009information,ding2015proof}.  A typical Spin glass instance can be achieved using random interactions graphs with $P(\bm J_{ij} = -1) = 0.5, P(\bm J_{ij} = 1) = 0.5$, and $\bm h_i = 0$.


\paragraph{Bijection of Ising and Boolean Optimization:}
It is valuable to observe that there is a bijection between Ising optimization (i.e., $\sigma \in \{-1,1\}$) and Boolean optimization (i.e., $x \in \{0,1\}$).  The transformation of $\sigma$-to-$x$ is given by:
\begin{subequations} \label{eq:spin2bool}
\begin{align}
    \sigma_i &= 2x_i - 1  ~\forall i \in {\cal N} \\ 
    \sigma_i\sigma_j &= 4x_ix_j - 2x_i - 2x_j + 1  ~\forall i,j \in {\cal E}
\end{align} \label{}
\end{subequations} 
and the inverse $x$-to-$\sigma$ is given by:
\begin{subequations}
\begin{align}
    x_i &= \frac{\sigma_i + 1}{2}  ~\forall i \in {\cal N} \\ 
    x_i x_j &= \frac{\sigma_i \sigma_j + \sigma_i + \sigma_j + 1}{4}  ~\forall i,j \in {\cal E} 
\end{align}
\end{subequations}
Consequently, any results from solving Ising models are also immediately applicable to the class of optimization problems referred to as Pseudo-Boolean Optimization or Quadratic Unconstrained Binary Optimization (QUBO):
%
\begin{align}
    & \min: \sum_{i,j \in {\cal E}} \bm c_{ij} x_i x_j + \sum_{i \in {\cal N}} \bm c_i x_i + \bm c  \label{eq:boolean_opt} \\
    & \mbox{s.t.: } x_i \in \{0, 1\}  ~\forall i \in {\cal N} \nonumber
\end{align}
%
%
In contrast to gate-based QC, which is Turing complete, QA specializes in optimizing Ising models.  The next section provides a brief introduction of how quantum mechanics are leveraged by QA to perform Ising model optimization.

\section{Foundations of Quantum Annealing}
\label{sec:qa_foundation}
Quantum annealing is an analog computing technique for minimizing discrete or continuous functions that takes advantage of the {\em exotic} properties of quantum systems.  This technique is particularly well-suited for finding optimal solutions of Ising models and has drawn significant interest due to hardware realizations via controllable quantum dynamical systems \cite{johnson2011quantum}.
Quantum annealing is composed of two key elements: leveraging quantum state to lift the minimization problem into an exponentially larger space, and slowly interpolating (i.e., annealing) between an initial easy problem and the target problem.  
The quantum lifting begins by introducing for each spin $\sigma_i \in\{-1,1\}$ a $2^N \times 2^N$ dimensional matrix $\widehat{\sigma}_i$ expressible as a Kronecker product of $N$ matrices of dimension $2\times 2$:
\begin{align}
    \widehat{\sigma}_i = 
    \underbrace{\begin{pmatrix}
    1 & 0 \\
    0 & 1
    \end{pmatrix}
        \mathop{\otimes}
    \cdots
        \mathop{\otimes}
    \begin{pmatrix}
    1 & 0 \\
    0 & 1
    \end{pmatrix}}_\text{$1$ to $i-1$}
        \mathop{\otimes}
    \underbrace{\begin{pmatrix}
    1 & 0 \\
    0 & -1
    \end{pmatrix}}_\text{$i^{\rm{th}}$ term} 
        \mathop{\otimes}
    \underbrace{\begin{pmatrix}
    1 & 0 \\
    0 & 1
    \end{pmatrix} 
        \mathop{\otimes}
    \cdots
        \mathop{\otimes}
    \begin{pmatrix}
    1 & 0 \\
    0 & 1
    \end{pmatrix}}_\text{$i+1$ to $N$}
\end{align}
In this lifted representation, the value of a spin $\sigma_i$ is identified with the two possible eigenvalues $1$ and $-1$ of the matrix $\widehat{\sigma}_i$. The quantum counterpart of the energy function defined in \eqref{eq:ising_eng} is the $2^N \times 2^N$ matrix obtained by substituting spins with the $\widehat{\sigma}$ matrices in the algebraic expression of the energy:
\begin{align}
    & \widehat{E} = \sum_{i,j \in {\cal E}} \bm J_{ij} \widehat{\sigma}_i \widehat{\sigma}_j + \sum_{i \in {\cal N}} \bm h_i \widehat{\sigma}_i  \label{eq:quantum_ising}
\end{align}
Notice that the eigenvalues of the matrix in \eqref{eq:quantum_ising} are the $2^N$ possible energy values obtained by evaluating the energy $E(\sigma)$ from $\eqref{eq:ising_eng}$ for all possible configurations of spins. This implies that finding the lowest eigenvalue of $\widehat{E}$ is tantamount to solving the minimization problem in \eqref{eq:ising_opt}.
This lifting is clearly impractical from the classical computing context as it transforms a minimization problem over $2^N$ configurations into computing the minimum eigenvalue of a $2^N \times 2^N$ matrix.
The key motivation for this approach is that it is possible to construct quantum systems with only $N$ quantum bits that attempt to find the minimum eigenvalue of this matrix. 

The annealing process provides a way of steering a quantum system into the a priori unknown eigenvector that minimizes the energy of \eqref{eq:quantum_ising} \cite{PhysRevE.58.5355,quant-ph-0001106}.  The core idea is to initialize the quantum system at the minimal eigenvector of a simple energy matrix $\widehat{E}_0$, for which an explicit formula is known. 
After the system is initialized, the energy matrix is interpolated from the easy problem to the target problem slowly over time.  Specifically, the energy matrix at a point during the anneal is given by $\widehat{E}_a(\Gamma) = (1-\Gamma)\widehat{E}_0 + \Gamma \widehat{E}$, with $\Gamma$ varying from $0$ to $1$.  When the anneal is complete, $\Gamma=1$ and the interactions in the quantum system are described by the target energy matrix.  The annealing time is the physical time taken by the system to evolve from $\Gamma=0$ to $\Gamma=1$. For suitable starting energy matrices $\widehat{E}_0$ and a sufficiently slow annealing time, theoretical results have demonstrated that a quantum system continuously remains at the minimal eigenvector of the interpolating matrix $\widehat{E}_a(\Gamma)$ \cite{albash2018adiabatic} and therefore achieves the minimum energy (i.e., a global optima) of the target problem.
Realizing this optimality result in practice has proven difficult due to corruption of the quantum system from the external environment.  Nevertheless, quantum annealing can serve as a heuristic for finding high-quality solutions to the Ising models, i.e.,   \eqref{eq:ising_opt}.


\subsection{Quantum Annealing Hardware}
\label{sec:qa_hardware}

\begin{figure}[t]
    \begin{center}
    \includegraphics[scale=0.90]{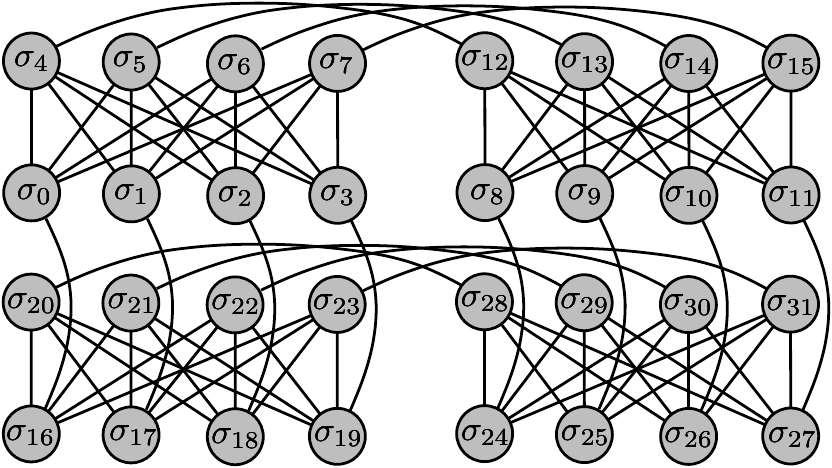}
    \end{center}
    \caption{A 2-by-2 Chimera graph illustrating the variable product limitations of D-Wave's 2000Q processor}
    \label{fig:chimera}
\end{figure}

Interest in the QA model is due in large part to D-Wave Systems, which has developed the first commercially available QA hardware platform \cite{johnson2011quantum}.  Given the computational challenges of classically simulating QA, this novel computing device represents the only viable method for studying QA at non-trivial scales, e.g., problems with more than 1000 qubits \cite{Boixo2014,PhysRevX.6.031015}.
At the most basic level, the D-Wave platform allows the user to program an Ising model by providing the parameters $\bm J, \bm h$ in \eqref{eq:ising_eng} and returns a collection of variable assignments from multiple annealing runs, which reflect optimal or near-optimal solutions to the input problem.

This seemingly simple interface is, however, hindered by a variety of constraints imposed by D-Wave's 2000Q hardware implementation. 
The most notable hardware restriction is the {\em 
Chimera} connectivity graph depicted in Figure \ref{fig:chimera}, where each edge indicates if the hardware supports a coupling term $\bm J_{ij}$ between a pair of qubits $i$ and $j$.  This sparse graph is a stark contrast to traditional quadratic optimization tools, where it is assumed that every pair of variables can interact.
%

The second notable hardware restriction is a limited coefficient programming range.  On the D-Wave 2000Q platform the parameters are constrained within the continuous parameter ranges of $-1 \leq \bm J_{ij} \leq 1$ and $-2 \leq \bm h_{i} \leq 2$.  At first glance these ranges may not appear to be problematic because the energy function \eqref{eq:ising_eng} can be rescaled into the hardware's operating range without any loss of generality.  However, operational realities of analog computing devices make the parameter values critically important to the overall performance of the hardware. These challenges include: persistent coefficient biases, which are an artifact of hardware slowly drifting out of calibration between re-calibration cycles; programming biases, which introduce some minor errors in the $\bm J, \bm h$ values that were requested; and environmental noise, which disrupts the quantum behavior of the hardware and results in a reduction of solution quality.  Overall, these hardware constraints have made the identification of QA-based performance gains notoriously challenging \cite{McGeoch:2013:EEA:2482767.2482797,PhysRevA.94.022337,10.1007/978-3-030-19212-9_11,1604.00319,ISAKOV2015265}.

Despite the practical challenges in using D-Wave's hardware platform, extensive experiments have suggested that QA can outperform some established local search methods (e.g., simulated annealing) on carefully designed Ising models \cite{1701.04579,PhysRevX.8.031016,PhysRevX.6.031015}.  However, demonstrating an unquestionable computational advantage over state-of-the-art methods on contrived and practical problems remains an open challenge.

\section{Methods for Ising Model Optimization}
\label{sec:algs}

The focus of this work is to compare and contrast the behavior of QA to a broad range of established optimization algorithms.  To that end, this work considers three core algorithmic categories: (1) {\em complete search} methods from the mathematical programming community; (2) {\em local search} methods developed by the statistical physics community; and (3) {\em quantum annealing} as realized by D-Wave's hardware platform.  
The comparison includes both state-of-the-art solution methods from the D-Wave benchmarking literature (e.g., Hamze-Freitas-Selby  \cite{HFS_impl_2017}, Integer Linear Programming  \cite{10.1007/978-3-030-19212-9_11}) and simple straw-man approaches (e.g., Greedy, Glauber Dynamics \cite{glauber1963time}, Min-Sum \cite{fossorier1999reduced,mezard2009information}) to highlight the solution quality of minimalist optimization approaches.  This section provides high-level descriptions of the algorithms; implementation details are available as open-source software \cite{ising_solvers,HFS_impl_2017}.

\subsection{Complete Search}

Unconstrained Boolean optimization, as in \eqref{eq:boolean_opt}, has been the subject of mathematical programming research for several decades \cite{BOROS2002155,Billionnet2007}.
This work considers the two most canonical formulations based on Integer Quadratic Programming and Integer Linear Programming.

\paragraph{Integer Quadratic Programming (IQP):}
This formulation consists of using black-box commercial optimization tools to solve \eqref{eq:boolean_opt} directly.  This model was leveraged in some of the first QA benchmarking studies \cite{McGeoch:2013:EEA:2482767.2482797} and received some criticism \cite{ibm_blog}.  However, the results presented here suggest that this model has become more competitive due to the steady progress of commercial optimization solvers.

\paragraph{Integer Linear Programming (ILP):}
This formulation is a slight variation of the IQP model where the variable products $x_ix_j$ are lifted into a new variable $x_{ij}$ and constraints are added to capture the conjunction $x_{ij} = x_i \wedge x_j$ as follows:
\begin{subequations}
\begin{align}
    & \min:  \sum_{i,j \in {\cal E}} \bm c_{ij} x_{ij} + \sum_{i \in {\cal N}} \bm c_i x_i + \bm c \\
    & \mbox{s.t.: } \nonumber \\
    & x_{ij} \geq x_{i} + x_{j} - 1, ~x_{ij} \leq x_{i}, ~x_{ij} \leq x_{j} ~\forall i,j \in {\cal E} \\
    & x_i \in \{0, 1\} ~\forall i \in {\cal N}, ~x_{ij} \in \{0, 1\} ~\forall i,j \in {\cal E} \nonumber 
\end{align}
\end{subequations}
This formulation was also leveraged in some of the first QA benchmarking studies \cite{ibm_blog,1306.1202} and \cite{Billionnet2007}, which suggest this is the best formulation for sparse graphs, as is the case with the D-Wave Chimera graph.  However, this work indicates that IQP solvers have improved sufficiently and this conclusion should be revisited.

%
%

\subsection{Local Search}

Although complete search algorithms are helpful in the validation of QA hardware \cite{baccari2018verification,10.1007/978-3-030-19212-9_11}, it is broadly accepted that local search algorithms are the most appropriate point of computational comparison to QA methods \cite{aaronson_blog2}.  Given that a comprehensive enumeration of local search methods would be a monumental undertaking, this work focuses on representatives from four distinct algorithmic categories including greedy, message passing, Markov Chain Monte Carlo, and large neighborhood search.

\paragraph{Greedy (GRD):}
The first heuristic algorithm considered by this work is a Steepest Coordinate Decent (SCD) greedy initialization approach.  This algorithm assigns the variables one-by-one, always taking the assignment that minimizes the objective value.  Specifically, the SCD approach begins with unassigned values, i.e., $\sigma_i = 0 ~\forall i \in {\cal N}$, and then repeatedly applies the following assignment rule until all of the variables have been assigned a value of $-1$ or 1:
\begin{subequations}
\begin{align}
    i, v &= \argmin_{i \in {\cal N}, v \in \{-1, 1\}} E(\sigma_1, \ldots, \sigma_{i-1}, v, \sigma_{i+1}, \ldots,\sigma_N) \\
    \sigma_i &= v
\end{align}
\end{subequations}
In each application, ties in the argmin are broken at random, giving rise to a potentially stochastic outcome of the heuristic.  Once all of the variables have been assigned, the algorithm is repeated until a runtime limit is reached and only the best solution found is returned.
Although this approach is very simple, it can be effective in Ising models with minimal amounts of frustration.

\paragraph{Message Passing (MP):}
The second algorithm considered by this work is a message-based Min-Sum (MS) algorithm \cite{fossorier1999reduced,mezard2009information}, which is an adaptation of the celebrated Belief Propagation algorithm for solving minimization problems on networks.
A key property of the MS approach is its ability to identify the global minimum of cost functions with a tree dependency structure between the variables; i.e., if no cycles are formed by the interactions in $\mathcal{E}$.  In the more general case of loopy dependency structures \cite{mezard2009information}, MS provides a heuristic minimization method.
It is nevertheless a popular technique favored in communication systems for its low computational cost and notable performance on random tree-like networks \cite{vuffray2014cavity}.

For the optimization model considered here, as in \eqref{eq:ising_opt}, the MS messages, $\epsilon_{i \rightarrow j}$,  are computed iteratively along \emph{directed} edges $i \rightarrow j$ and $j \rightarrow i$ for each edge $(i,j)\in \mathcal{E}$, according to the Min-Sum equations:
\begin{subequations}
\begin{align}
    \epsilon^{t+1}_{i \rightarrow j} = \mathrm{SSL}(2\bm J_{ij},2\bm h_i + \sum_{k \in  \mathcal{E}(i) \setminus j}\epsilon^{t}_{k \rightarrow j} ) \label{eq:min_sum)} \\
    \mathrm{SSL}(x,y) = \min(x,y)-\min(-x,y) -x
\end{align}
\end{subequations}
Here, $\mathcal{E}(i) \setminus j$ denotes the neighbors of $i$ without $j$ and $\mathrm{SSL}$ denotes the Symmetric Saturated Linear transfer function.
Once a fix-point of \eqref{eq:min_sum)} is obtained or a prescribed runtime limit is reached, the MS algorithm outputs a configuration based on the following formula:
\begin{align}
    \sigma_{i} = - \mathrm{sign}\left( 2\bm h_i + \sum_{k \in \mathcal{E}(i)}\epsilon_{k \rightarrow j} \right) \label{eq:min_sum_assignement}
\end{align}
By convention, if the argument of the {\em sign} function is 0, a value of $1$ or $-1$ is assigned randomly with equal probability.

\paragraph{Markov Chain Monte Carlo (MCMC):}
MCMC algorithms include a wide range of methods to generate samples from complex probability distributions. A natural Markov Chain for the Ising model is given by Glauber dynamics, where the value of each variable is updated according to its conditional probability distribution. Glauber dynamics is often used as a method for producing samples from Ising models at {\em finite temperature} \cite{glauber1963time}. This work considers the so-called {\em Zero Temperature} Glauber Dynamics (GD) algorithm, which is the optimization variant of the Glauber dynamics sampling method, and which is also used in physics as a simple model for describing avalanche phenomena in magnetic materials \cite{dhar1997zero}. From the optimization perspective, this approach is a 
single-variable greedy local search algorithm.

A step $t$ of the GD algorithm consists in checking each variable $i\in \mathcal{N}$ in a random order and comparing the objective cost of the current configuration $\sigmab^t$ to the configuration with the variable $\sigma^{t}_i$ being flipped. If the objective value is lower in the flipped configuration, i.e.,
$E(\sigmab^t) > E(\sigma^{t}_1,\ldots,-\sigma^{t}_i,\ldots,\sigma^{t}_N)$,
then the flipped configuration is selected as the new current configuration $\sigmab^{t+1} = (\sigma^{t}_1,\ldots,-\sigma^{t}_i,\ldots,\sigma^{t}_N)$.  When the objective difference is 0, the previous or new configuration is selected randomly with equal probability. 
If after visiting all of the variables, no one single-variable flip can improve the current assignment, then the configuration is identified as a local minimum and the algorithm is restarted with a new randomly generated configuration.  This process is repeated until a runtime limit is reached.

\paragraph{Large Neighborhood Search (LNS):}
The state-of-the-art meta-heuristic for benchmarking D-Wave-based QA algorithms is the Hamze-Freitas-Selby (HFS) algorithm \cite{Hamze:2004:FT:1036843.1036873,1409.3934}.  The core idea of this algorithm is to extract low treewidth subgraphs of the given Ising model and then use dynamic programming to quickly compute the optimal configuration of these subgraphs.  This extract and optimize process is repeated until a specified time limit is reached.  This approach has demonstrated remarkable results in a variety of benchmarking studies \cite{1604.00319,king2015performance,1701.04579,10.1007/978-3-030-19212-9_11,junger2019performance}.
The notable success of this solver can be attributed to three key factors.  First, it is highly specialized to solving Ising models on the Chimera graphs (i.e., Figure \ref{fig:chimera}), a topological structure that is particularly amenable to low treewidth subgraphs.  Second, it leverages integer arithmetic instead of floating point, which provides a significant performance improvement but also leads to notable precision limits.  Third, the baseline implementation is a highly optimized C code \cite{HFS_impl_2017}, which runs at near-ideal performance.

\subsection{Quantum Annealing}

Extending the theoretical overview from Section \ref{sec:qa_foundation}, the following implementation details are required to leverage the D-Wave 2000Q platform as a reliable optimization tool.
The QA algorithm considered here consists of programming the Ising model of interest and then repeating the annealing process some number of times (i.e., {\em num\_reads}) and then returning the lowest energy solution that was found among all of those replicates.  No correction or solution polishing is applied in this solver.  By varying the number of reads considered (e.g., from 10 to 10,000), the solution quality and total runtime of the QA algorithm increases.  It is important to highlight that the D-Wave platform provides a wide variety of parameters to control the annealing process (e.g., annealing time, qubit offsets, custom annealing schedules, etc.).  In the interest of simplicity and reproducibility, this work does not leverage any of those advanced features and it is likely that the results presented here would be further improved by careful utilization of those additional capabilities \cite{PhysRevA.96.042322,Adame_2020,PhysRevApplied.11.044083}.

Note that all of the problems considered in this work have been generated to meet the implementation requirements discussed in Section \ref{sec:qa_hardware} for a specific D-Wave chip deployed at Los Alamos National Laboratory.  Consequently, no problem transformations are required to run the instances on the target hardware platform.  Most notably, no embedding or rescaling is required.
This approach is standard practice in QA evaluation studies and the arguments for it are discussed at length in \cite{coffrin2016challenges,10.1007/978-3-030-19212-9_11}.

\section{Structure Detection Experiments}
\label{sec:structure}

This section presents the primary result of this work.  Specifically, it analyzes three crafted optimization problems of increasing complexity—the Biased Ferromagnet, Frustrated Biased Ferromagnet, and Corrupted Biased Ferromagnet—all of which highlight the potential for QA to quickly identify the global structural properties of these problems.  The algorithm performance analysis focuses on two key metrics, solution quality over time (i.e., performance profile) and the minimum hamming distance to any optimal solution over time.  The hamming distance metric is particularly informative in this study as the  problems have been designed to have local minima that are very close to the global optimum in terms of objective value, but are very distant in terms of hamming distance.  The core finding is that QA produces solutions that are close to global optimality, both in terms of objective value and hamming distance.

\paragraph{Problem Generation:}
All problems considered in this work are defined by simple probabilistic graphical models and are generated on a specific D-Wave hardware graph. To avoid bias towards one particular random instance, 100 instances are generated and the mean over this collection of instances is presented.  Additionally, a random gauge transformation is applied to every instance to obfuscate the optimal solution and mitigate artifacts from the choice of initial condition in each solution approach.

\paragraph{Computation Environment:}
The CPU-based algorithms are run on HPE ProLiant XL170r servers with dual Intel 2.10GHz CPUs and 128GB memory.  Gurobi 9.0 \cite{gurobi} was used for solving the Integer Programming (ILP/IQP) formulations.  All of the algorithms were configured to only leverage one thread and the reported runtime reflects the wall clock time of each solver's core routine and does not include pre-processing or post-processing of the problem data.

The QA computation is conducted on a D-Wave 2000Q quantum annealer deployed at Los Alamos National Laboratory.  This computer has a 16-by-16 Chimera cell topology with random omissions; in total, it has 2032 spins (i.e., {$\cal N$}) and 5924 couplers (i.e., {$\cal E$}).  The hardware is configured to execute 10 to 10,000 annealing runs using a 5-microsecond annealing time per run and a random gauge transformation every 100 runs, to mitigate the various sources of bias in the problem encoding.  The reported runtime of the QA hardware reflects the amount of {\em on-chip time} used; it does not include the overhead of communication or scheduling of the computation, which takes about one to two seconds.  Given a sufficient engineering effort to reduce overheads, on-chip time would be the dominating runtime factor.

\begin{figure}[t]
    \begin{center}
    \includegraphics[width=0.46\textwidth]{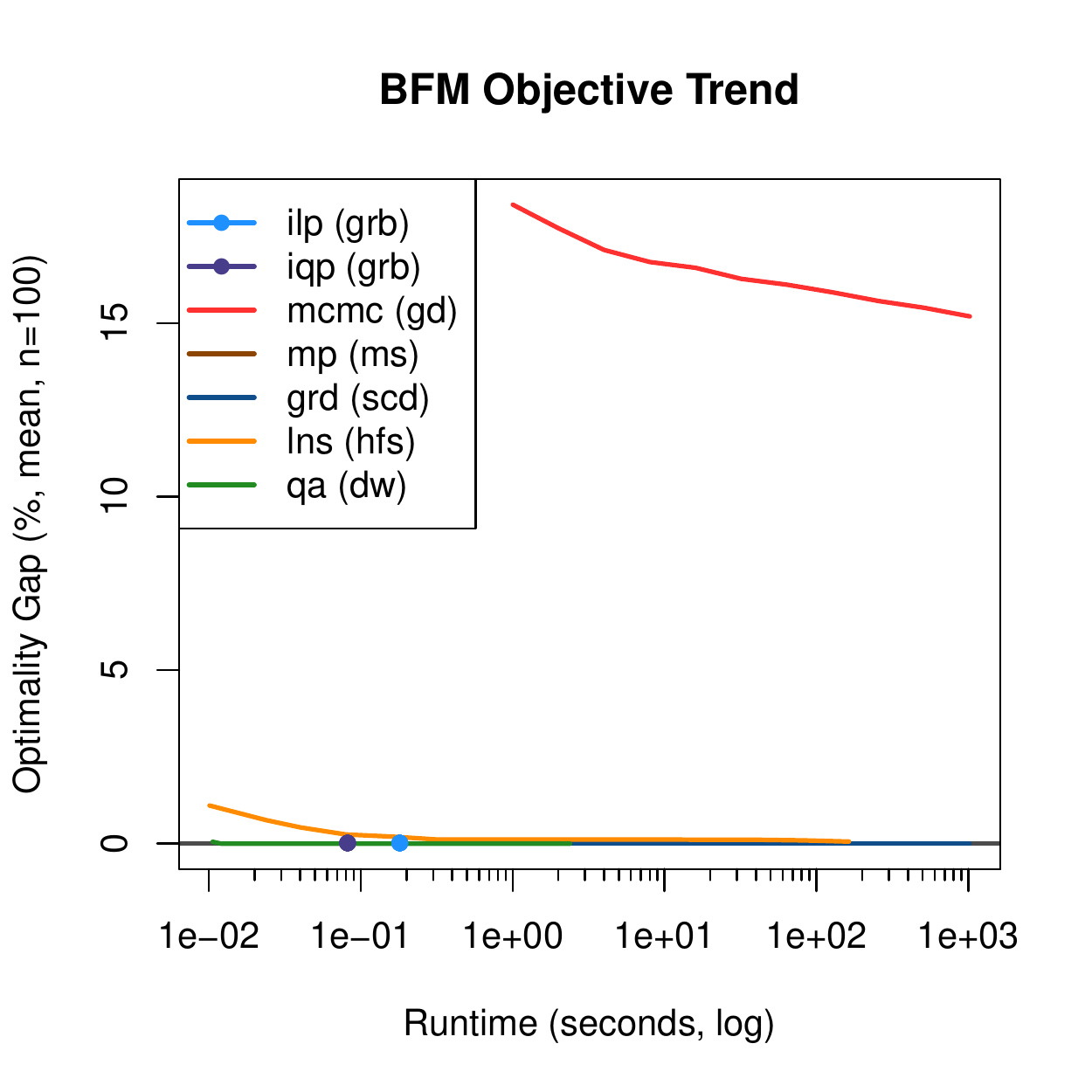}
    \includegraphics[width=0.46\textwidth]{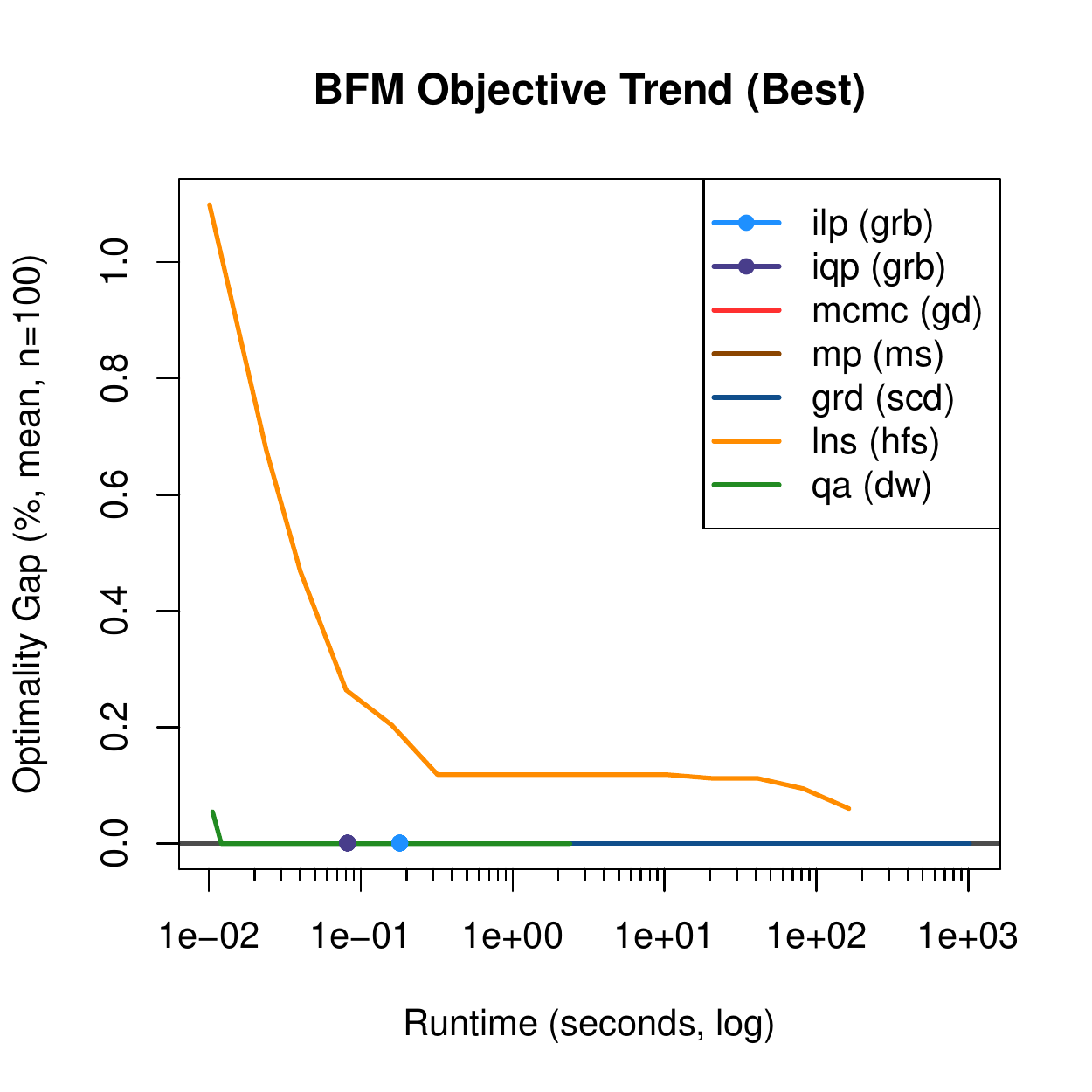}\\
    \vspace{-0.3cm}
    \includegraphics[width=0.46\textwidth]{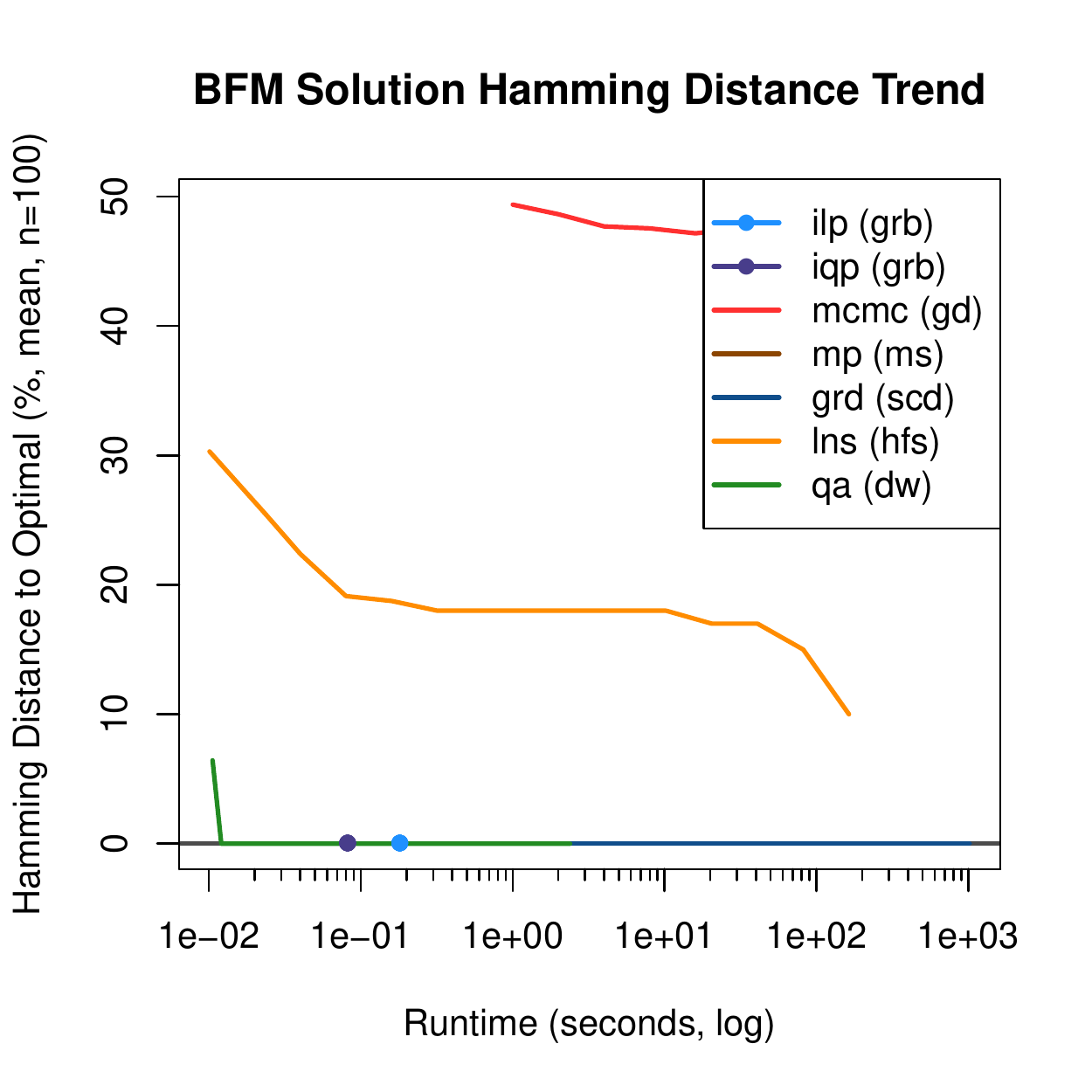}
    \includegraphics[width=0.46\textwidth]{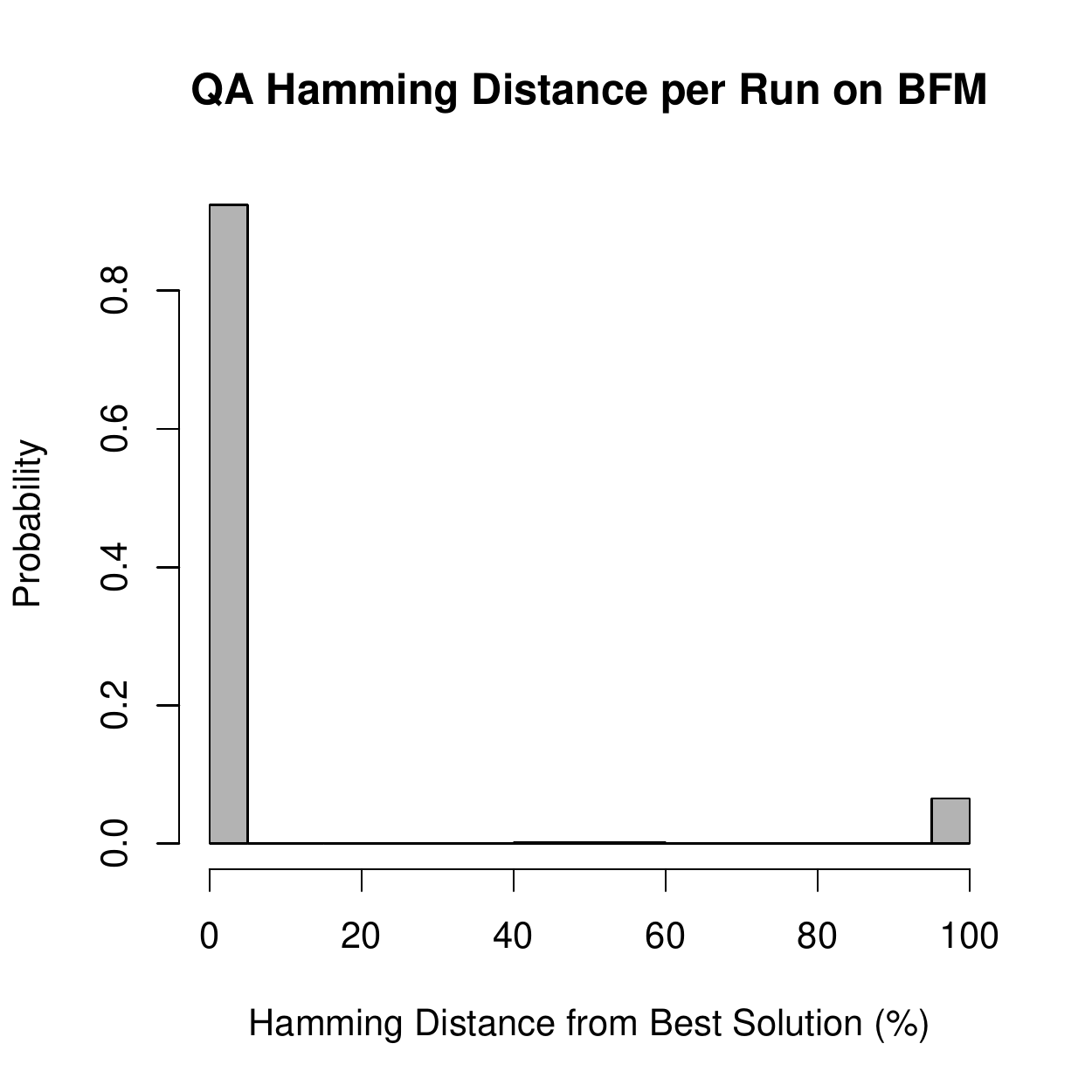}
    \end{center}
    \vspace{-0.5cm}
    \caption{Performance profile (top) and Hamming Distance (bottom) analysis for the Biased Ferromagnet instance}
    \label{fig:bfm}
\end{figure}

\subsection{The Biased Ferromagnet}
\vspace{-0.2cm}
\begin{align}
    \label{eq:bfm} \tag{BFM}
    \bm J_{ij} &= -1.00 ~\forall i,j \in {\cal E};\\ 
    P(\bm h_i = 0.00) = 0.990 &, P(\bm h_i = -1.00) = 0.010  ~\forall i \in {\cal N} \nonumber
\end{align}

Inspired by the Ferromagnet model, this study begins with Biased FerroMagnet \eqref{eq:bfm} model—a toy problem to build an intuition for a type of structure that QA can exploit.  Notice that this model has no frustration and has a few linear terms that bias it to prefer $\sigma_i = 1$ as the global optimal solution.  W.h.p. $\sigma_i = 1$ is a unique optimal solution and the assignment of $\sigma_i = -1$  is a local minimum that is sub-optimal by $0.02 \cdot |{\cal N}|$ in expectation and has a maximal hamming distance of $|{\cal N}|$.  The local minimum is an attractive solution because it is nearly optimal; however, it is hard for a local search solver to escape from it due to its hamming distance from the true global minimum.  This instance presents two key algorithmic challenges: first, one must effectively detect the global structure (i.e., all the variables should take the same value); second, one must correctly discriminate between the two nearly optimal solutions that are very distant from one another.

Figure \ref{fig:bfm} presents the results of running all of the algorithms from Section \ref{sec:algs} on the BFM model.  The key observations are as follows:
\begin{itemize}
    \item Both the greedy (i.e., SCD) and relaxation-based solvers (i.e., IQP/ILP/MS) correctly identify this problem's structure and quickly converge on the globally optimal solution (Figure \ref{fig:bfm}, top-right).
    \item Neighborhood-based local search methods (e.g., GD) tend to get stuck in the local minimum of this problem.  Even advanced local search methods (e.g., HFS) may miss the global optimum in rare cases (Figure \ref{fig:bfm}, top).
    \item The hamming distance analysis indicates that QA has a high probability (i.e., greater than 0.9) of finding the exact global optimal solution (Figure \ref{fig:bfm}, bottom-right).  This explains why just 20 runs is sufficient for QA to find the optimal solution w.h.p. (Figure \ref{fig:bfm}, top-right).
\end{itemize}
A key observation from this toy problem is that making a continuous relaxation of the problem (e.g., IQP/ILP/MS) can help algorithms detect global structure and avoid local minima that present challenges for neighborhood-based local search methods (e.g., GD/LNS).  QA has comparable performance to these relaxation-based methods, both in terms of solution quality and runtime, and does appear to detect the global structure of the BFM problem class.

However encouraging these results are, the BFM problem is a straw-man that is trivial for five of the seven solution methods considered here.  The next experiment introduces frustration to the BFM problem to understand how that impacts problem difficulty for the solution methods considered.

\begin{figure}[t]
    \begin{center}
    \includegraphics[width=0.46\textwidth]{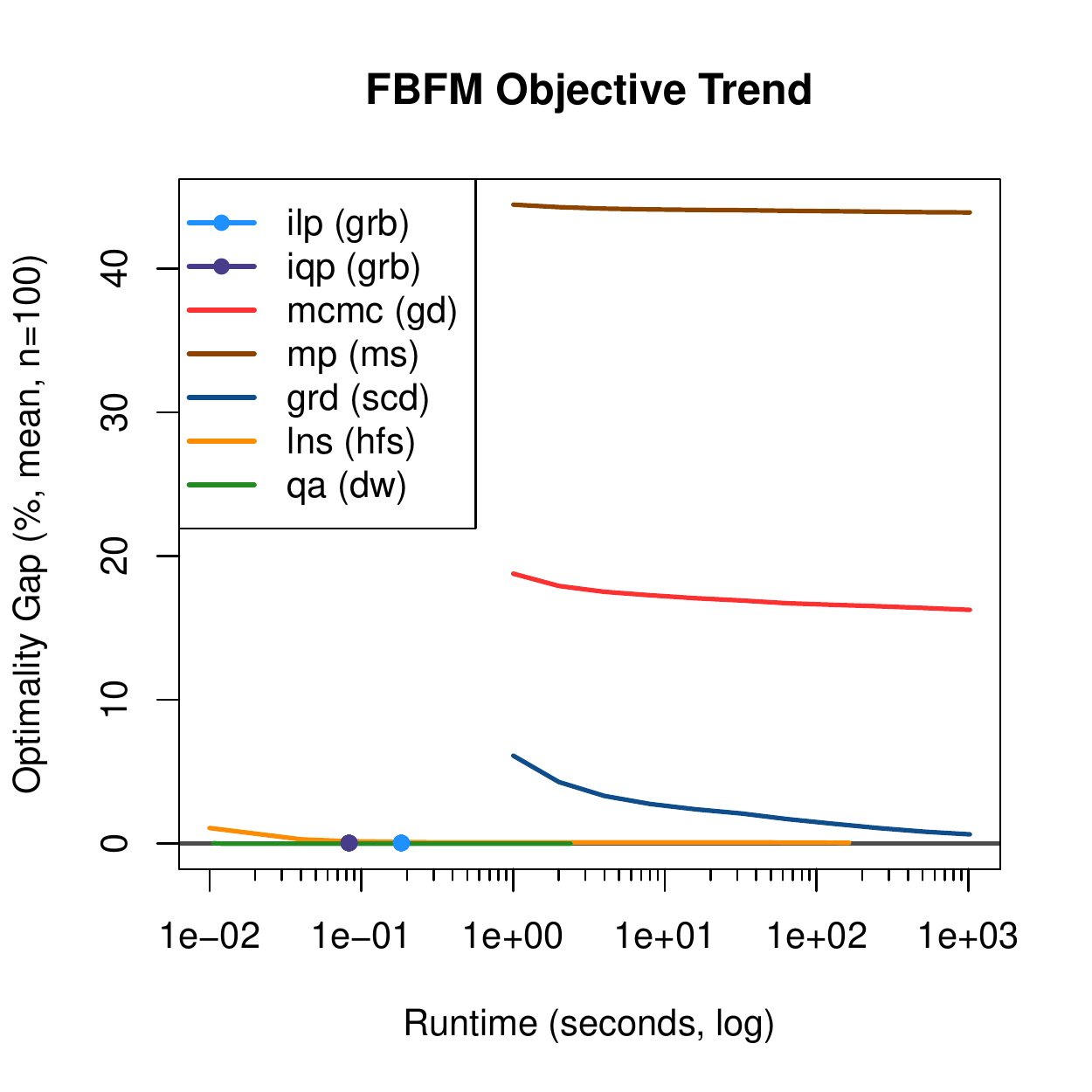}
    \includegraphics[width=0.46\textwidth]{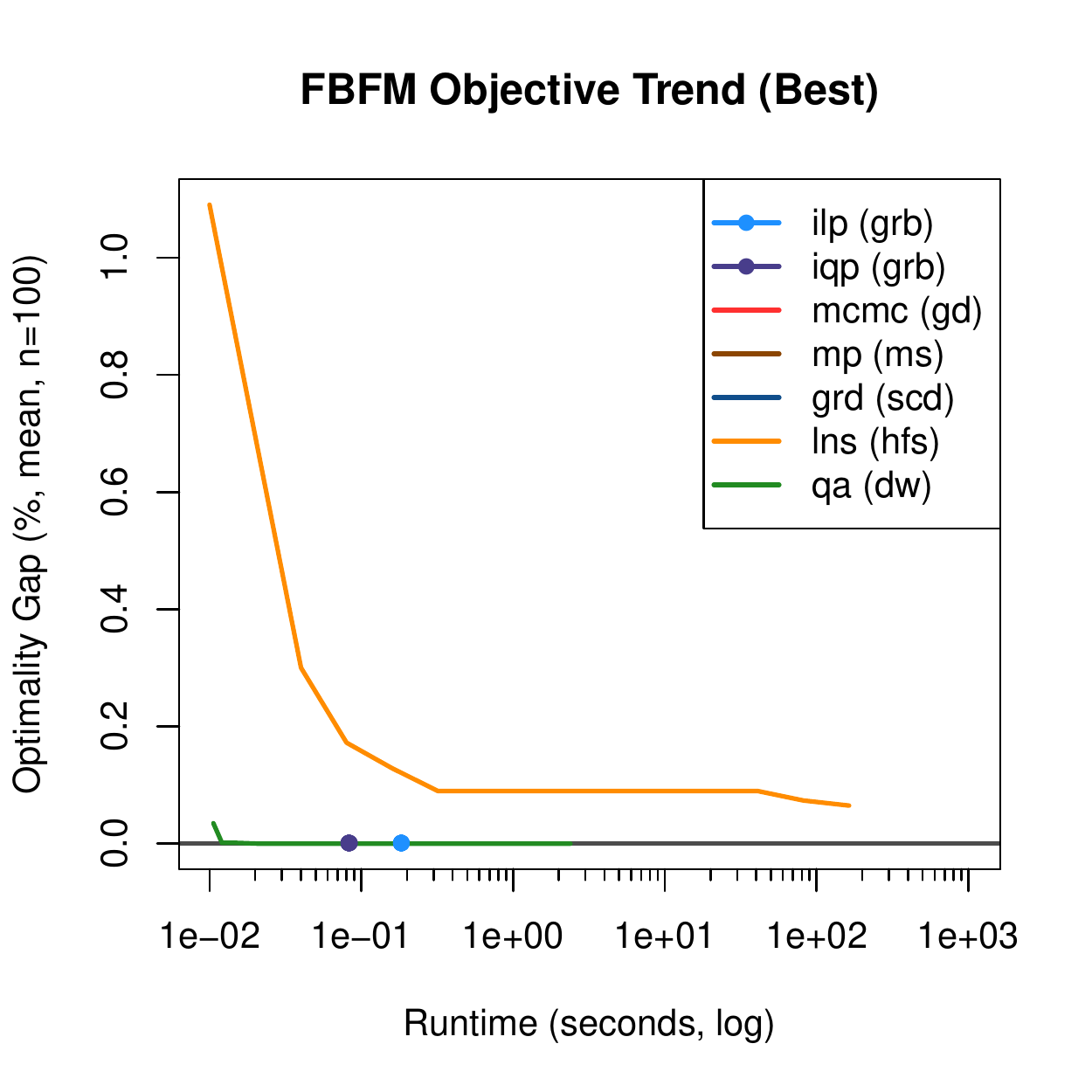}\\
    \vspace{-0.3cm}
    \includegraphics[width=0.46\textwidth]{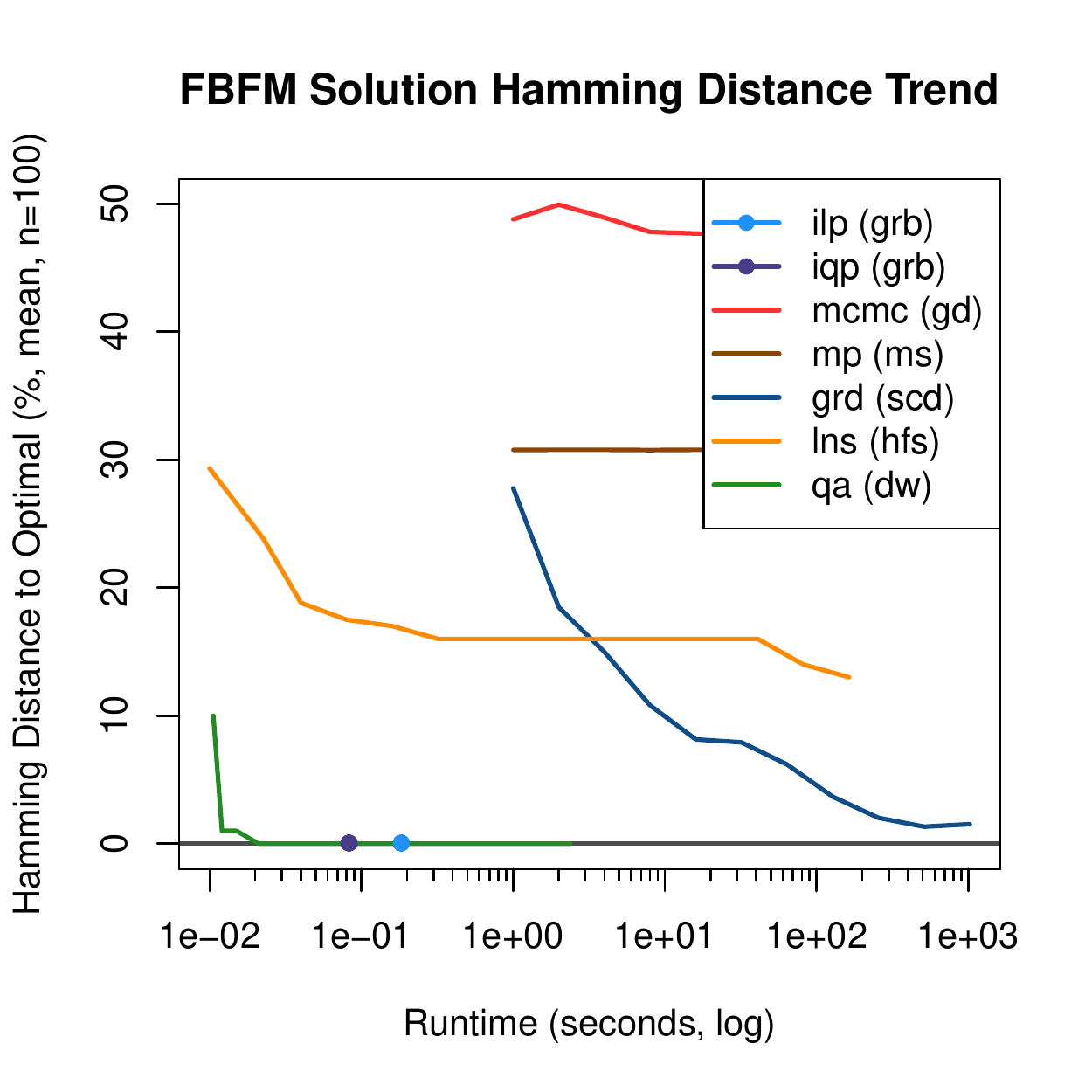}
    \includegraphics[width=0.46\textwidth]{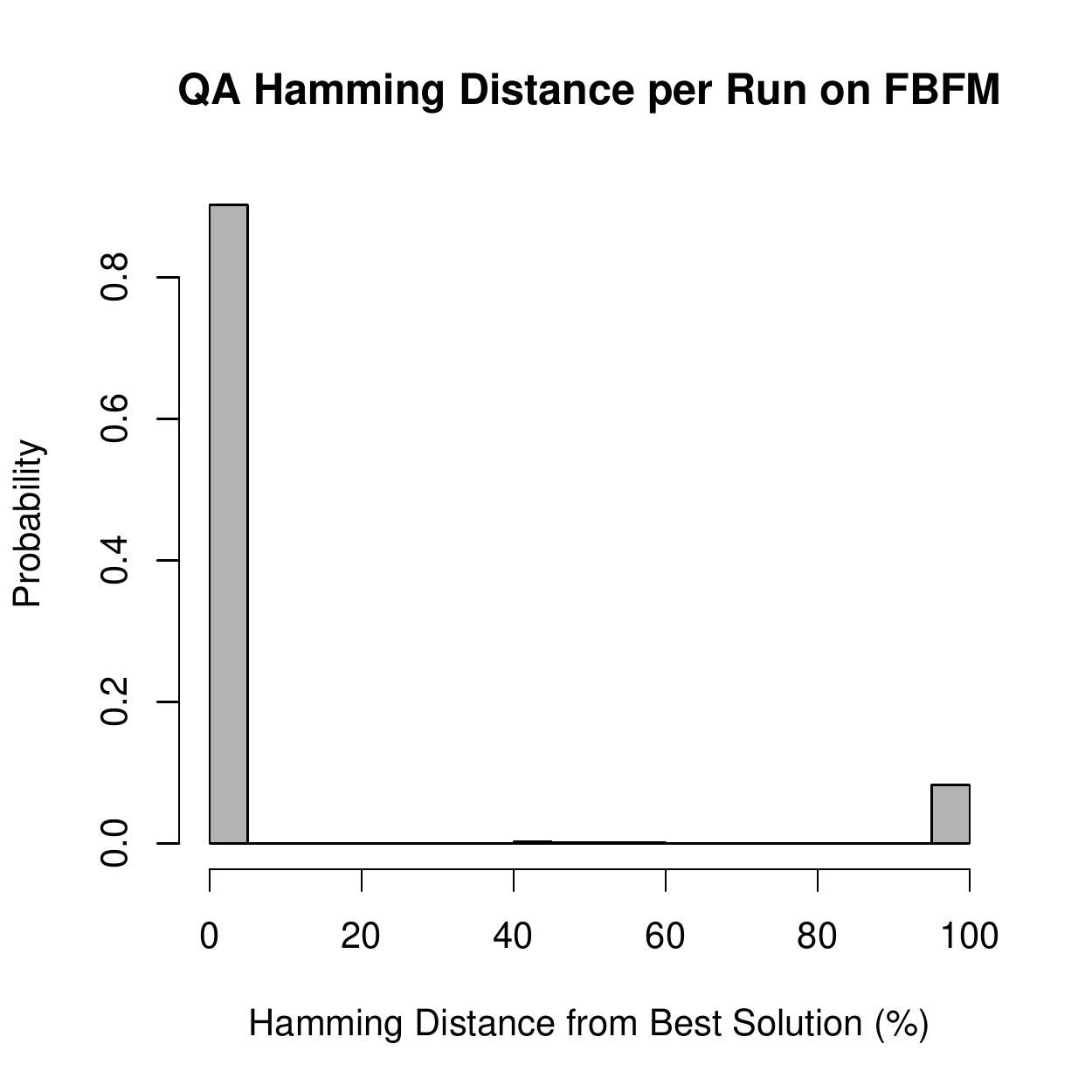}
    \end{center}
    \vspace{-0.5cm}
    \caption{Performance profile (top) and Hamming Distance (bottom) analysis for the Frustrated Biased Ferromagnet instance}
    \label{fig:fbfm}
\end{figure}

\subsection{The Frustrated Biased Ferromagnet}
\vspace{-0.2cm}
\begin{align}
    \label{eq:fbfm} \tag{FBFM}
    \bm J_{ij} &= -1.00 ~\forall i,j \in {\cal E} \\
    P(\bm h_i = 0.00) = 0.970, P(\bm h_i = & -1.00) = 0.020, P(\bm h_i = 1.00) = 0.010 ~\forall i \in {\cal N} \nonumber
\end{align}

The next step considers a slightly more challenging problem called a Frustrated Biased Ferromagnet \eqref{eq:fbfm}, which is a specific case of the random field Ising model \cite{d1985random} and similar in spirit to the Clause Problems considered in \cite{dwave_ocp}.  The FBFM deviates from the BFM by introducing frustration among the linear terms of the problem.  Notice that on average 2\% of the decision variables locally prefer $\sigma_i = 1$ while 1\% prefer $\sigma_i = -1$.  Throughout the optimization process these two competing preferences must be resolved, leading to frustration.  W.h.p. this model has the same unique global optimal solution as the BFM that occurs when $\sigma_i = 1$.  The opposite assignment of $\sigma_i = -1$  remains a local minimum that is sub-optimal by $0.02 \cdot |{\cal N}|$ in expectation and has a maximal hamming distance of $|{\cal N}|$.  By design, the energy difference of these two extreme assignments is consistent with BFM, to keep the two problem classes as similar as possible.

Figure \ref{fig:fbfm} presents the same performance analysis for the FBFM model.  The key observations are as follows:
\begin{itemize}
    \item When compared to BFM, FBFM presents an increased challenge for the simple greedy (i.e., SCD) and local search (i.e., GD/MS) algorithms.
    \item Although the SCD algorithm is worse than HFS in terms of objective quality, it is comparable or better in terms of hamming distance (Figure \ref{fig:fbfm}, bottom-left).  This highlights how these two metrics capture different properties of the underlying algorithms.
    \item The results of QA and the relaxation-based solvers (i.e., IQP/ILP), are nearly identical to the BFM case, suggesting that this type of frustration does not present a significant challenge for these solution approaches.
\end{itemize}
These results suggest that frustration in the linear terms alone (i.e., $\bm h$) is not sufficient for building optimization tasks that are non-trivial for a wide variety of general purpose solution methods.  In the next study, frustration in the quadratic terms (i.e., $\bm J$) is incorporated to increase the difficulty for the relaxation-based solution methods.

\subsection{The Corrupted Biased Ferromagnet}
\vspace{-0.2cm}
\begin{align}
    \label{eq:cbfm} \tag{CBFM}
    P(\bm J_{ij} = -1.00) = 0.625, P(\bm J_{ij} = 0.20) &= 0.375 ~\forall i,j \in {\cal E} \\
    P(\bm h_i = 0.00) = 0.970, P(\bm h_i = -1.00) = 0.020 &, P(\bm h_i = 1.00) = 0.010 ~\forall i \in {\cal N} \nonumber
\end{align}

The inspiration for this instance is to leverage insights from the theory of Spin glasses to build more computationally challenging problems.  The core idea is to carefully corrupt the ferromagnetic problem structure  with frustrating anti-ferromagnetic links that obfuscate the ferromagnetic properties without completely destroying them.  A parameter sweep of different corruption values yields the Corrupted Biased FerroMagnet \eqref{eq:cbfm} model, which retains the global structure that $\sigma_i = 1$ is a near globally optimal solution w.h.p., while obfuscating this property with misleading  anti-ferromagnetic links and frustrated local fields.
\begin{figure}[t]
    \begin{center}
    \includegraphics[width=0.46\textwidth]{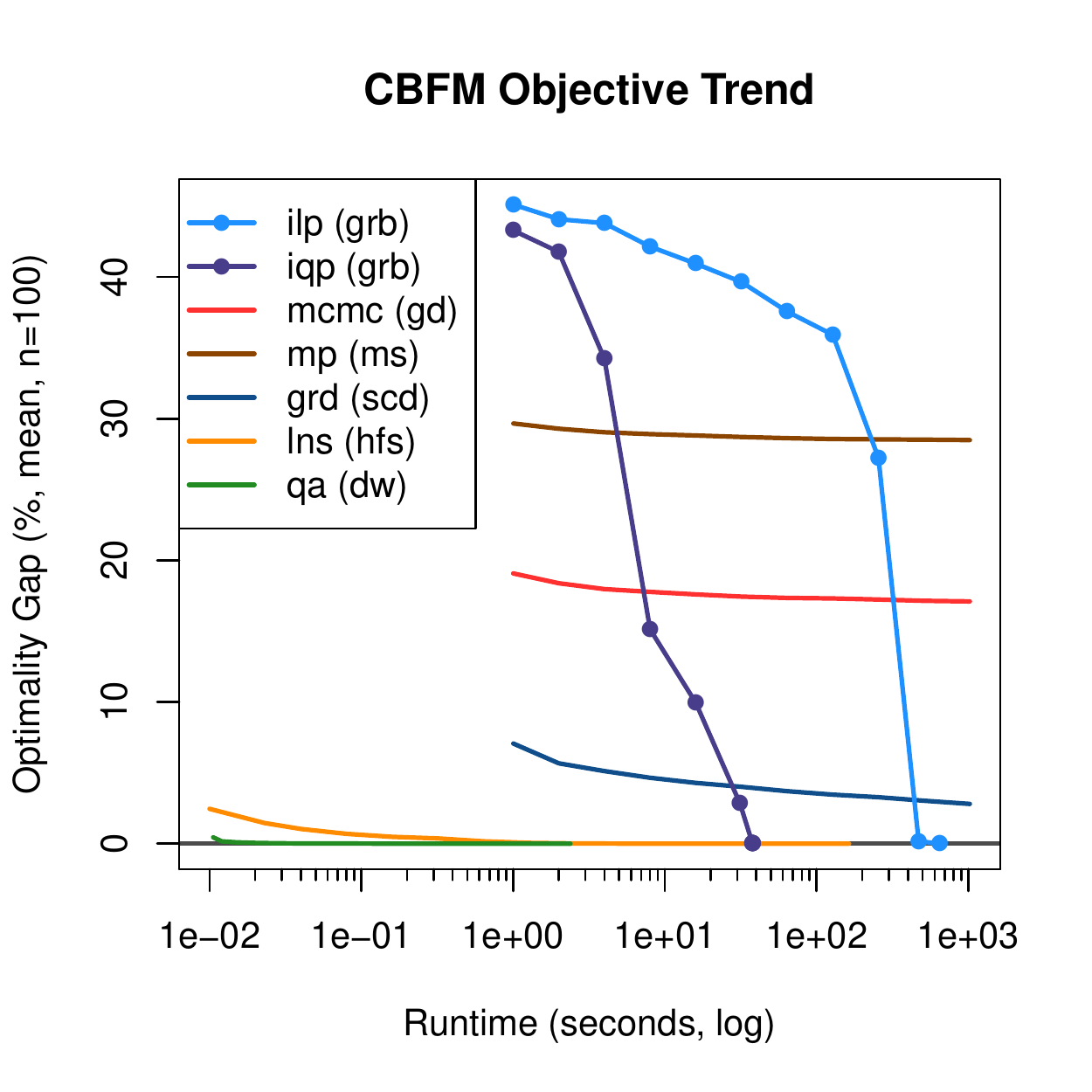}
    \includegraphics[width=0.46\textwidth]{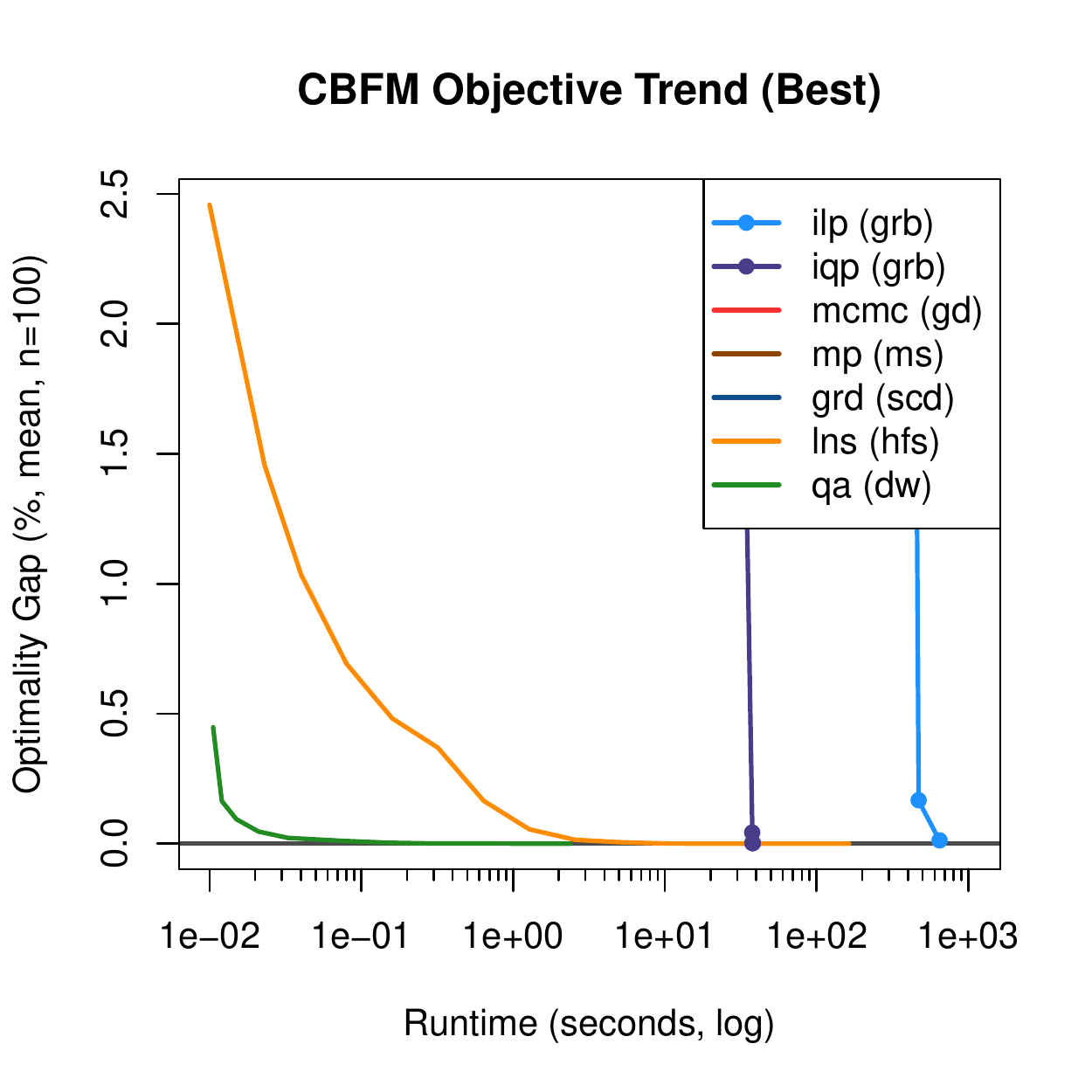}\\
    \vspace{-0.3cm}
    \includegraphics[width=0.46\textwidth]{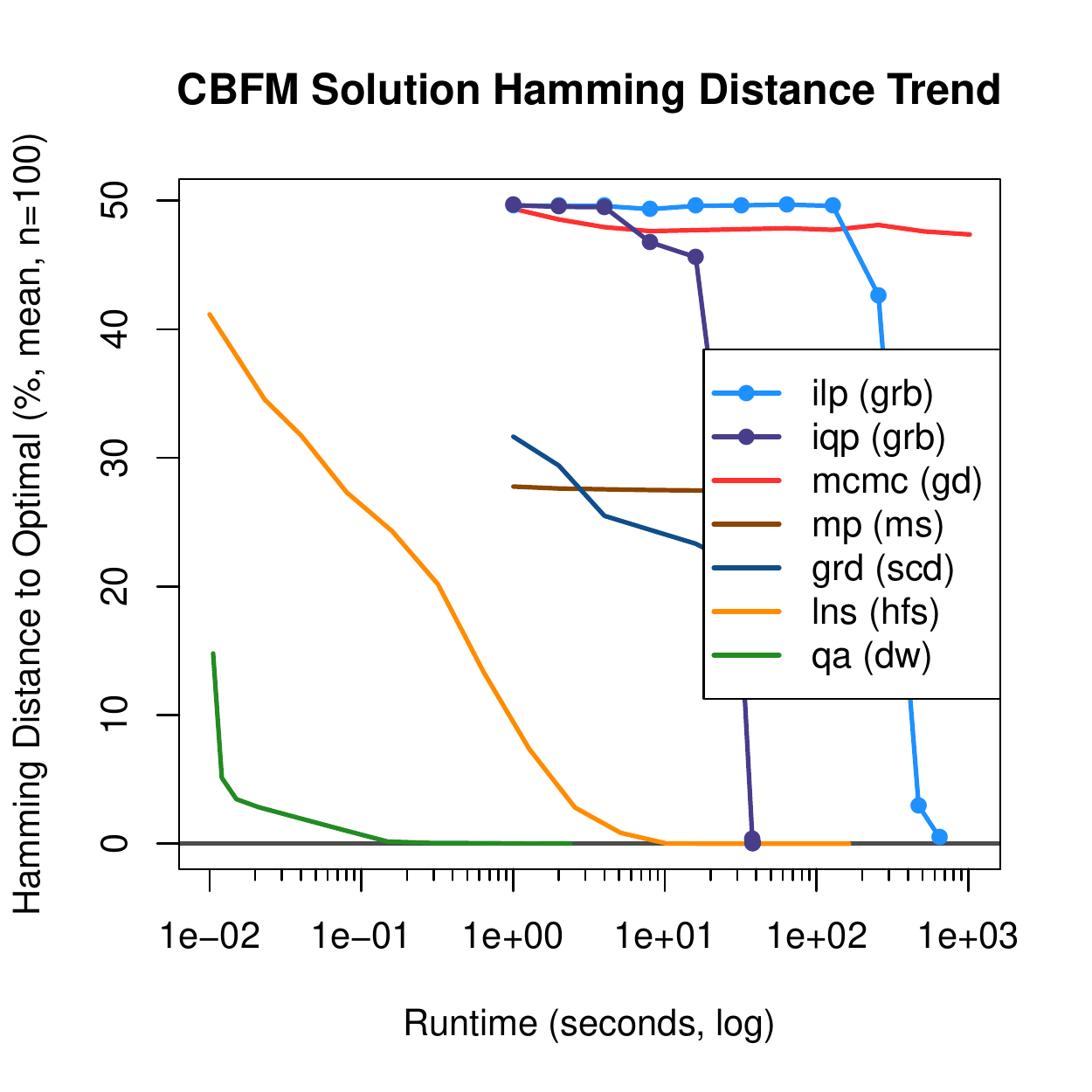}
    \includegraphics[width=0.46\textwidth]{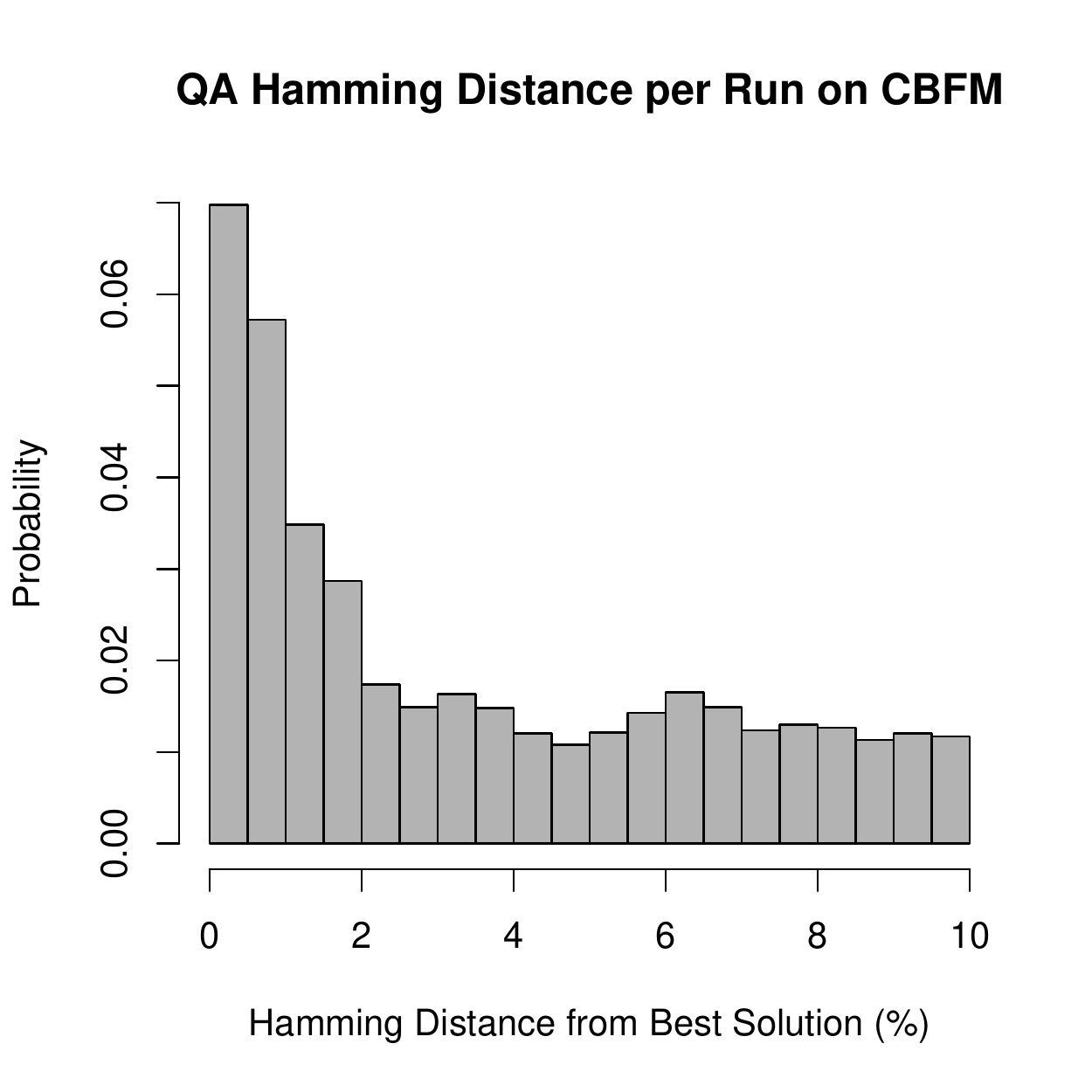}
    \end{center}
    \vspace{-0.5cm}
    \caption{Performance profile (top) and Hamming Distance (bottom) analysis for the Corrupted Biased Ferromagnet instance}
    \label{fig:cbfm}
\end{figure}
Figure \ref{fig:cbfm} presents a similar performance analysis for the CBFM model.  The key observations are as follows:
\begin{itemize}
    \item In contrast to the BFM and FBFM cases, solvers that leverage continuous relaxations, such as IQP and ILP, do not immediately identify this problem's structure and can take between 50 to 700 seconds to identify the globally optimal solution (Figure \ref{fig:cbfm}, top-left).
    \item The advanced local search method (i.e., HFS) consistently converges to a global optimum (Figure \ref{fig:cbfm}, top-right), which does not always occur in the BFM and FBFM cases.
    \item Although the MS algorithm is notably worse than GD in terms of objective quality, it is notably better in terms of hamming distance. This further indicates how these two metrics capture different properties of the underlying algorithms (Figure \ref{fig:cbfm}, bottom-left).
    \item Although this instance presents more of a challenge for QA than BFM and FBFM, QA still finds the global minimum with high probability; 500-1000 runs is sufficient to find a near-optimal solution in all cases.  This is 10 to 100 times faster than the next-best algorithm, HFS (Figure \ref{fig:cbfm}, top-right).
    \item The hamming distance analysis suggests that the success of the QA approach is that it has a significant probability (i.e., greater than $0.12$) of returning a solution that has a hamming distance of less than 1\% from the global optimal solution (Figure \ref{fig:cbfm}, bottom-right).
\end{itemize}
The overarching trend of this study is that QA is successful in detecting the global structure of the BFM, FBFM, and CBFM instances (i.e., low hamming distance to optimal, w.h.p.).  Furthermore, it can do so notably faster than all of the other algorithms considered here.  This suggests that, in this class of problems, QA brings a unique value that is not captured by the other algorithms considered.  Similar to how the relaxation methods succeed at the BFM and FBFM instances, we hypothesize that the success of QA on the CBFM instance is driven by the solution search occurring in a smooth high-dimensional continuous space as discussed in Section \ref{sec:qa_foundation}.  In this instance class, QA may also benefit from so-called {\em finite-range tunnelling} effects, which allows QA to change the state of multiple variables simultaneously (i.e., global moves) \cite{Farhi472,PhysRevX.6.031015}.
Regardless of the underlying cause, QA's performance on the CBFM instance is  particularly notable and worthy of further investigation.

\begin{figure}[t]
    \begin{center}
    \includegraphics[width=0.46\textwidth]{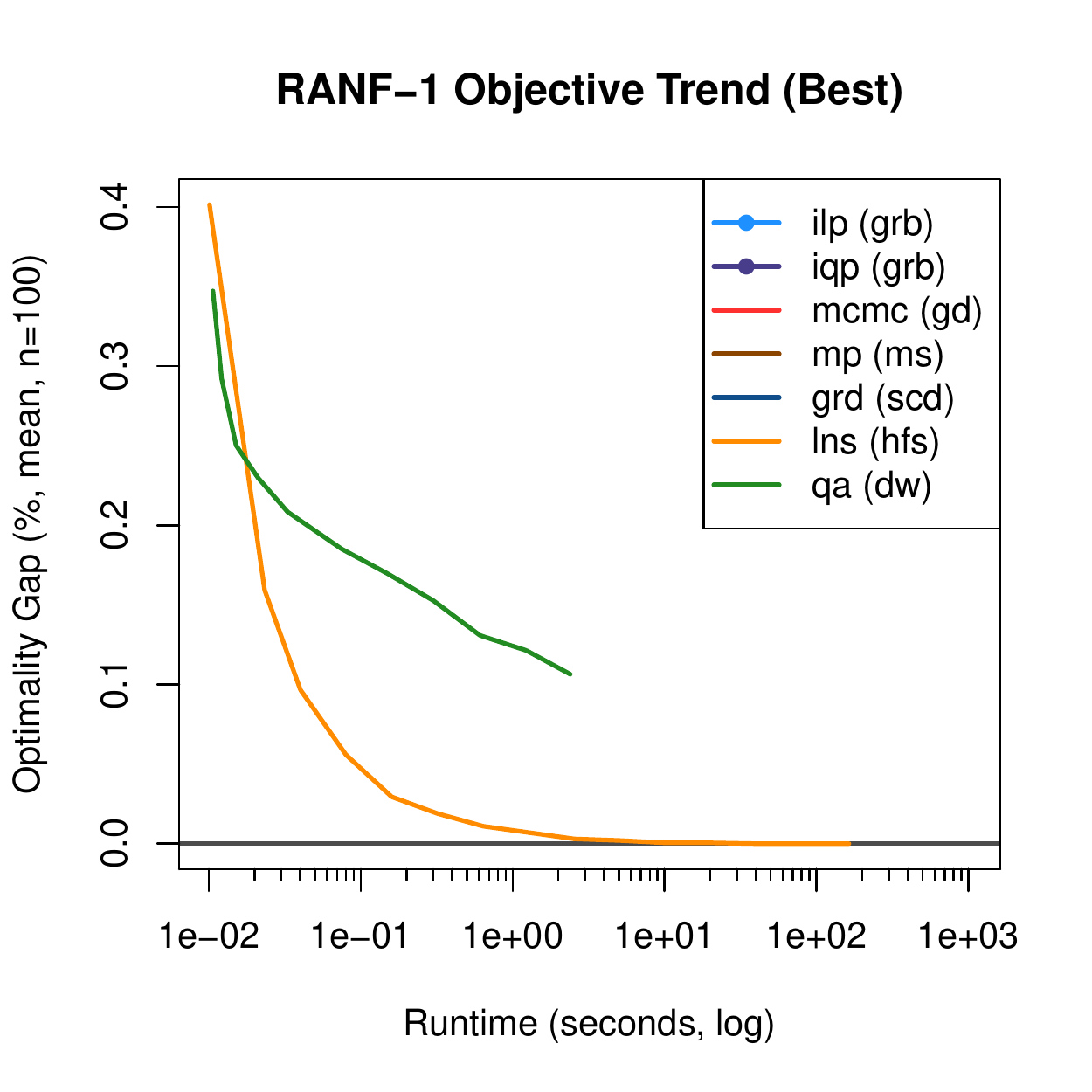}
    \includegraphics[width=0.46\textwidth]{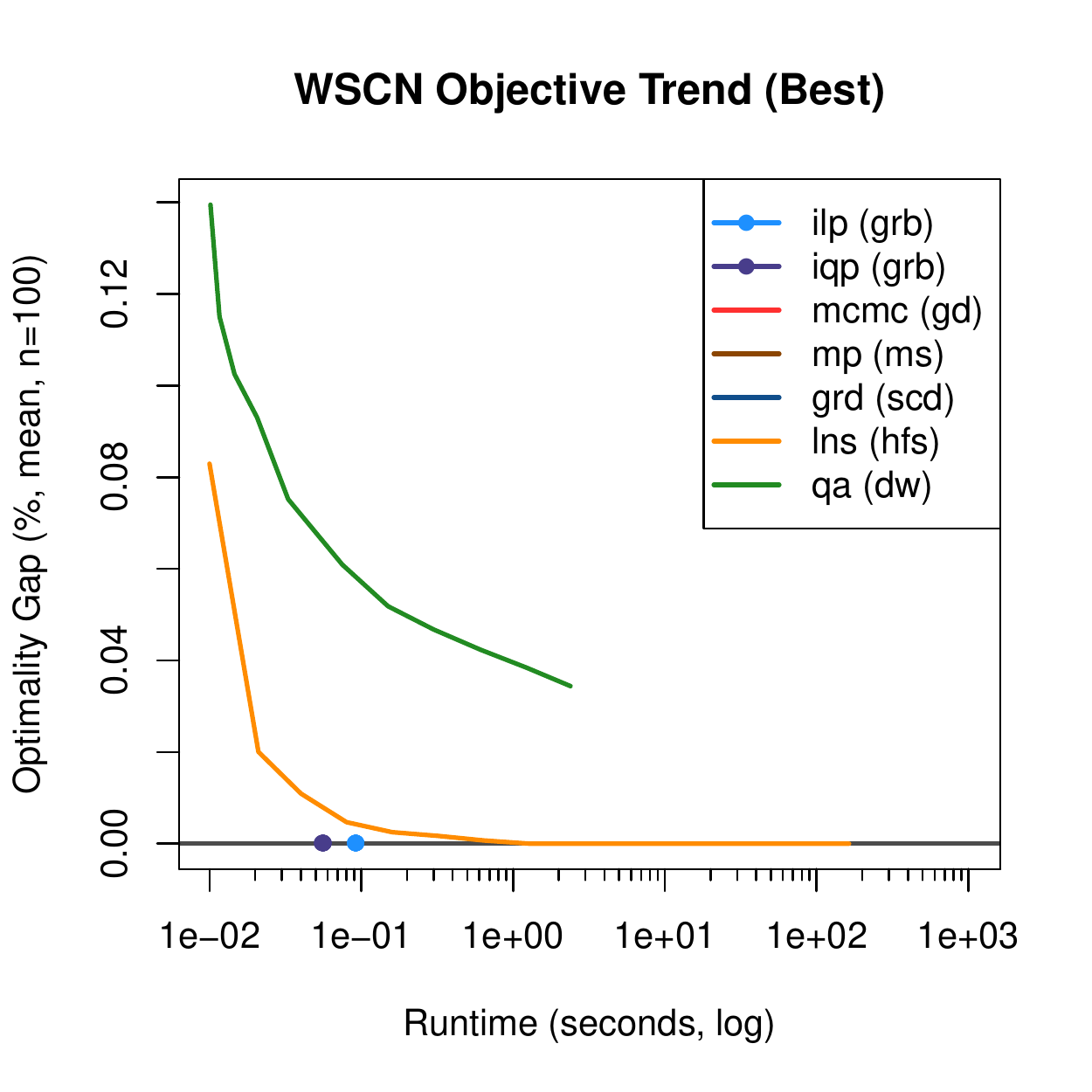}\\
    \vspace{-0.3cm}
    \includegraphics[width=0.46\textwidth]{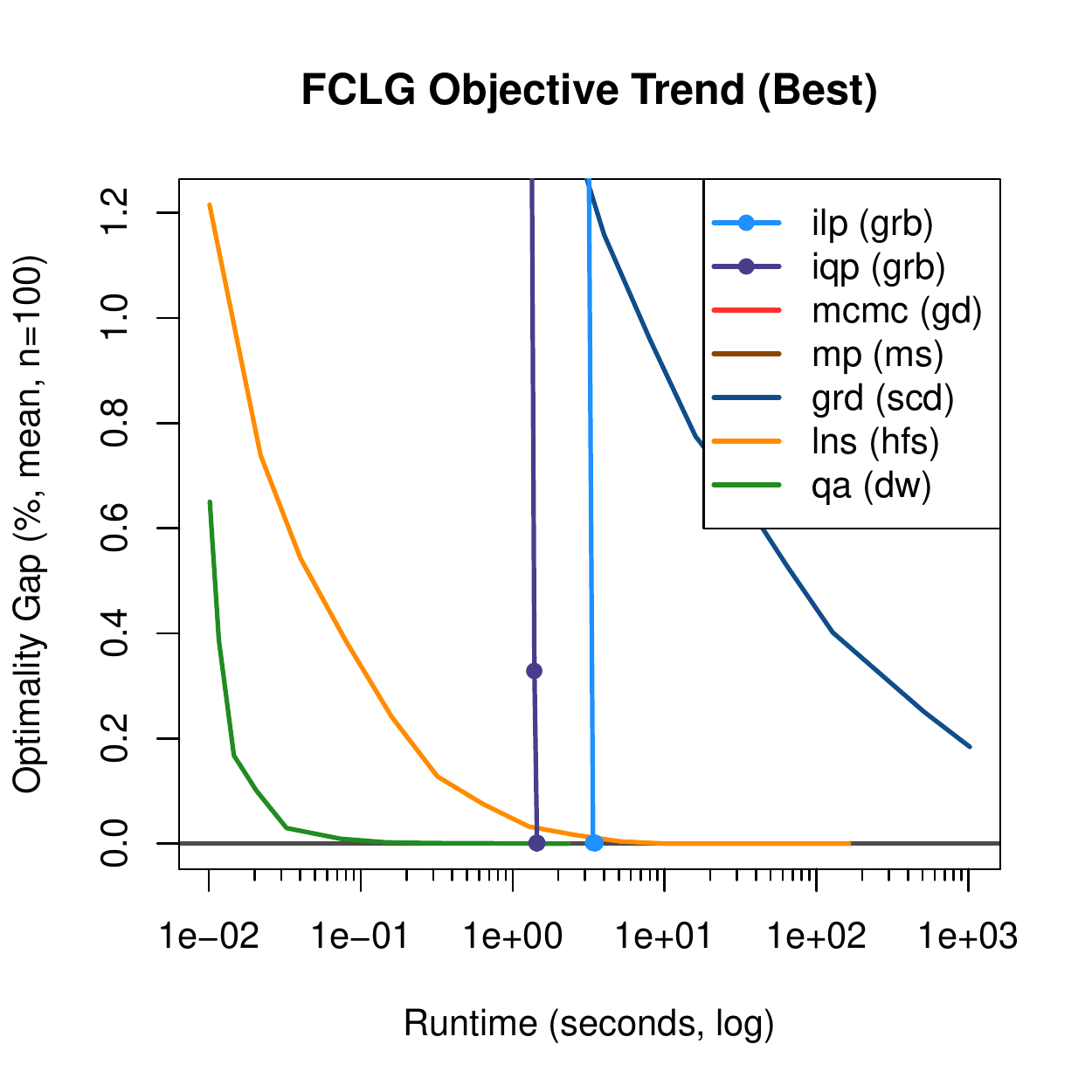}
    \end{center}
    \vspace{-0.5cm}
    \caption{Performance profiles of other problem classes from the literature}
    \label{fig:other-inst}
\end{figure}

\subsection{Bias Structure Variants}

As part of the design process uniform field variants of the problems proposed herein were also considered.  These variants featured weaker and more uniform distributed bias terms.  Specifically, the term $P(\bm h_i = -1.00) = 0.010$ was replaced with $P(\bm h_i = -0.01) = 1.000$.  Upon continued analysis, it was observed that the stronger and less-uniform bias terms resulted in more challenging cases for all of the solution methods considered, and hence, were selected as the preferred design for the problems proposed by this work.  In the interest of completeness, Appendix \ref{apx:bias} provides a detailed analysis of the uniform-field variants of the BFM, FBFM, and CBFM instances to illustrate how this problem variant impacts the performance of the solution methods considered here.

\subsection{A Comparison to Other Instance Classes}

The CBFM problem was designed to have specific structural properties that are beneficial to the QA approach.  It is important to note that not all instance classes have such an advantageous structure.  This point is highlighted in Figure \ref{fig:other-inst}, which compares three landmark problem classes from the QA benchmarking literature: Weak-Strong Cluster Networks (WSCN) \cite{PhysRevX.6.031015}, Frustrated Cluster Loops with Gadgets (FCLG) \cite{PhysRevX.8.031016}, and Random Couplers and Fields (RANF-1) \cite{10.1007/978-3-030-19212-9_11,1306.1202}.
These results show that D-Wave's current 2000Q hardware platform can be outperformed by local and complete search methods on some classes of problems.  However, it is valuable to observe that these previously proposed instance classes are either relatively easy for local search algorithms (i.e., WSCN and RANF) or relatively easy for complete search algorithms (i.e., WSCN and FCLG), both of which are not ideal properties for conducting benchmarking studies.  To the best of our knowledge, the proposed CBFM problem is the first instance class that presents a notable computational challenge for both local search and complete search algorithms.

\section{Quantum Annealing as a Primal Heuristic}
\label{sec:hybrid}

\begin{figure}[t]
    \begin{center}
    \includegraphics[width=0.46\textwidth]{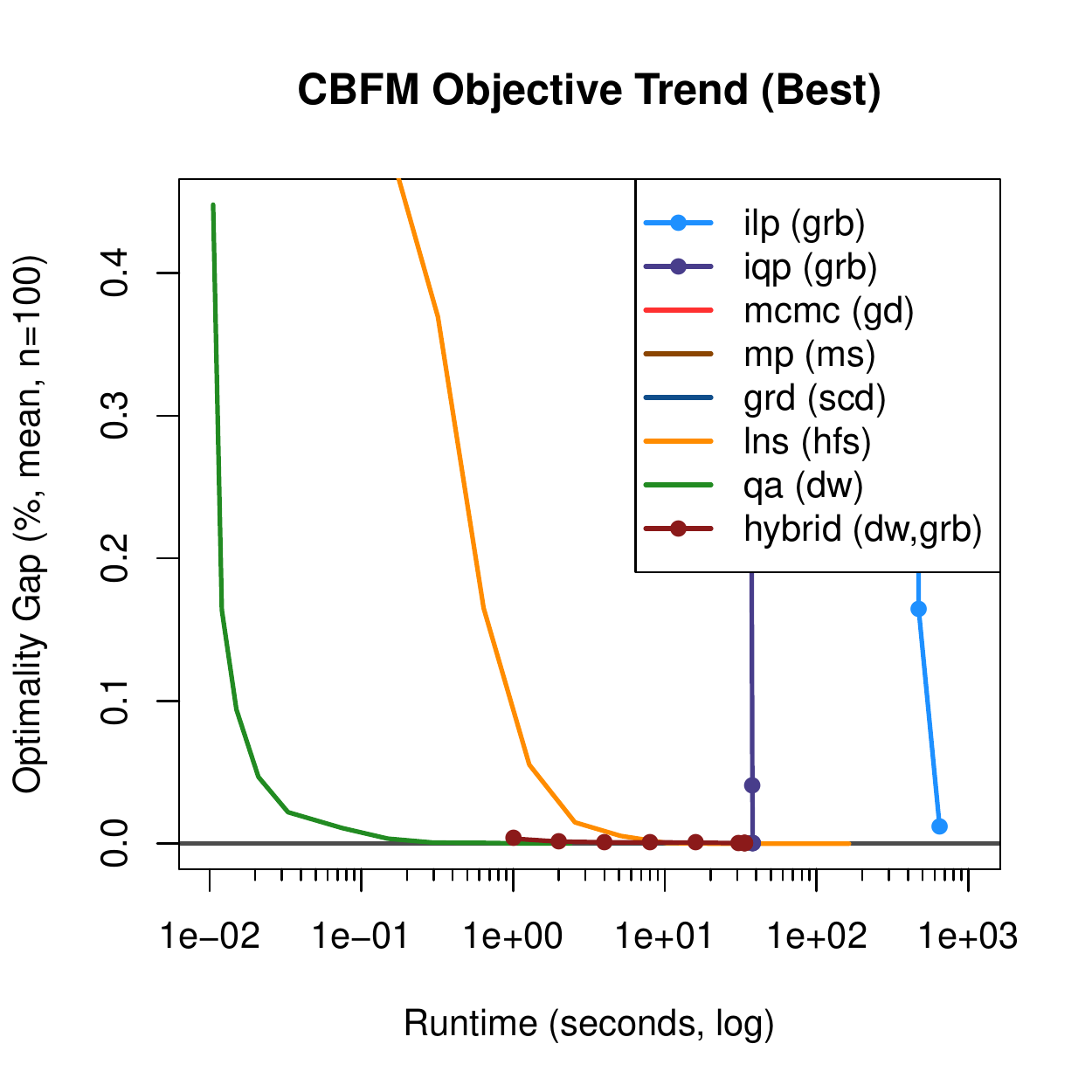}
    \end{center}
    \vspace{-0.5cm}
    \caption{Performance profile of Warm-Starting IQP with QA solutions}
    \label{fig:hybrid}
\end{figure}

QA's notable ability to find high-quality solutions to the CBFM problem suggests the development of hybrid algorithms, which leverage QA for finding upper bounds within a complete search method that can also provide global optimality proofs.  A simple version of such an approach was developed where 1000 runs of QA were used to warm-start the IQP solver with a high-quality initial solution.  The results of this hybrid approach are presented in Figure \ref{fig:hybrid}.  The IQP solver clearly benefits from the warm-start on short time scales.  However, it does not lead to a notable reduction in the time to producing the optimality proof.  This suggests that a state-of-the-art hybrid complete search solver needs to combine QA for finding upper bounds with more sophisticated lower-bounding techniques, such as those presented in \cite{baccari2018verification,junger2019performance}.

\section{Conclusion}
\label{sec:conclusion}

This work explored how quantum annealing hardware might be able to support heuristic algorithms in finding high-quality solutions to challenging combinatorial optimization problems.  A careful analysis of quantum annealing's performance on the Biased Ferromagnet, Frustrated Biased Ferromagnet, and Corrupted Biased Ferromagnet problems with more than 2,000 decision variables suggests that this approach is capable of quickly identifying the structure of the optimal solution to these problems, while a variety of local and complete search algorithms struggle to identify this structure.  This result suggests that integrating quantum annealing into meta-heuristic algorithms could yield unique variable assignments and increase the discovery of high-quality solutions.

Although demonstration of a runtime advantage was not the focus of this work, the success of quantum annealing on the Corrupted Biased Ferromagnet problem compared to other solution methods is a promising outcome for QA and warrants further investigation.  An in-depth theoretical study of the Corrupted Biased Ferromagnet case could provide deeper insights into the structural properties that quantum annealing is exploiting in this problem and would provide additional insights into the classes of problems that have the best chance to demonstrate an unquestionable computational advantage for quantum annealing hardware.
It is important to highlight that while the research community is currently searching for an unquestionable computational advantage for quantum annealing hardware by any means necessary, significant additional research will be required to bridge the gap between contrived hardware-specific optimization tasks and practical optimization applications.

\section*{Acknowledgments}
The research presented in this work was supported by the Laboratory Directed Research and Development program of Los Alamos National Laboratory under project numbers  20180719ER and 20190195ER.

\clearpage
\appendix

\section*{Appendix}

\begin{figure}[t]
    \begin{center}
    \includegraphics[width=0.43\textwidth]{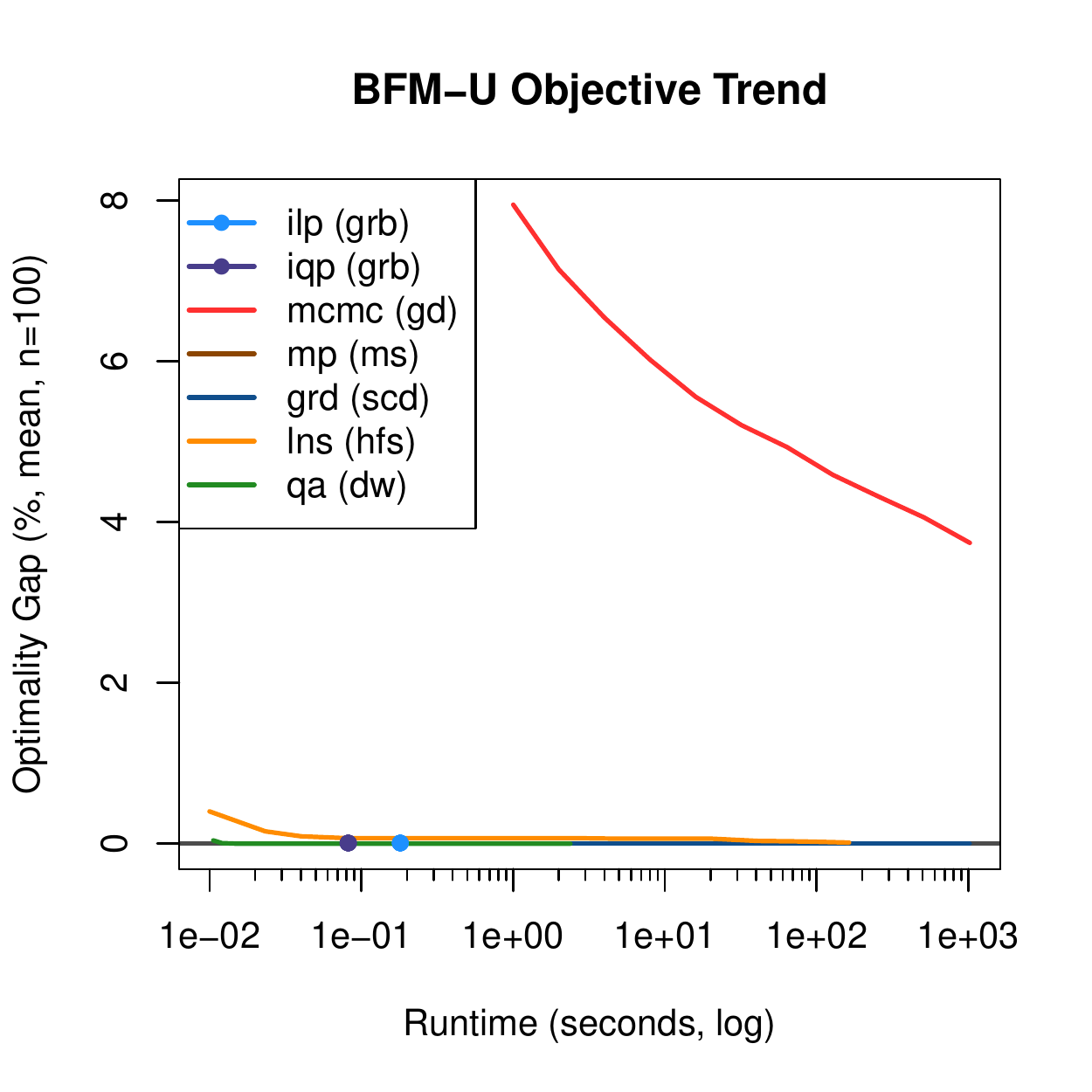}
    \includegraphics[width=0.43\textwidth]{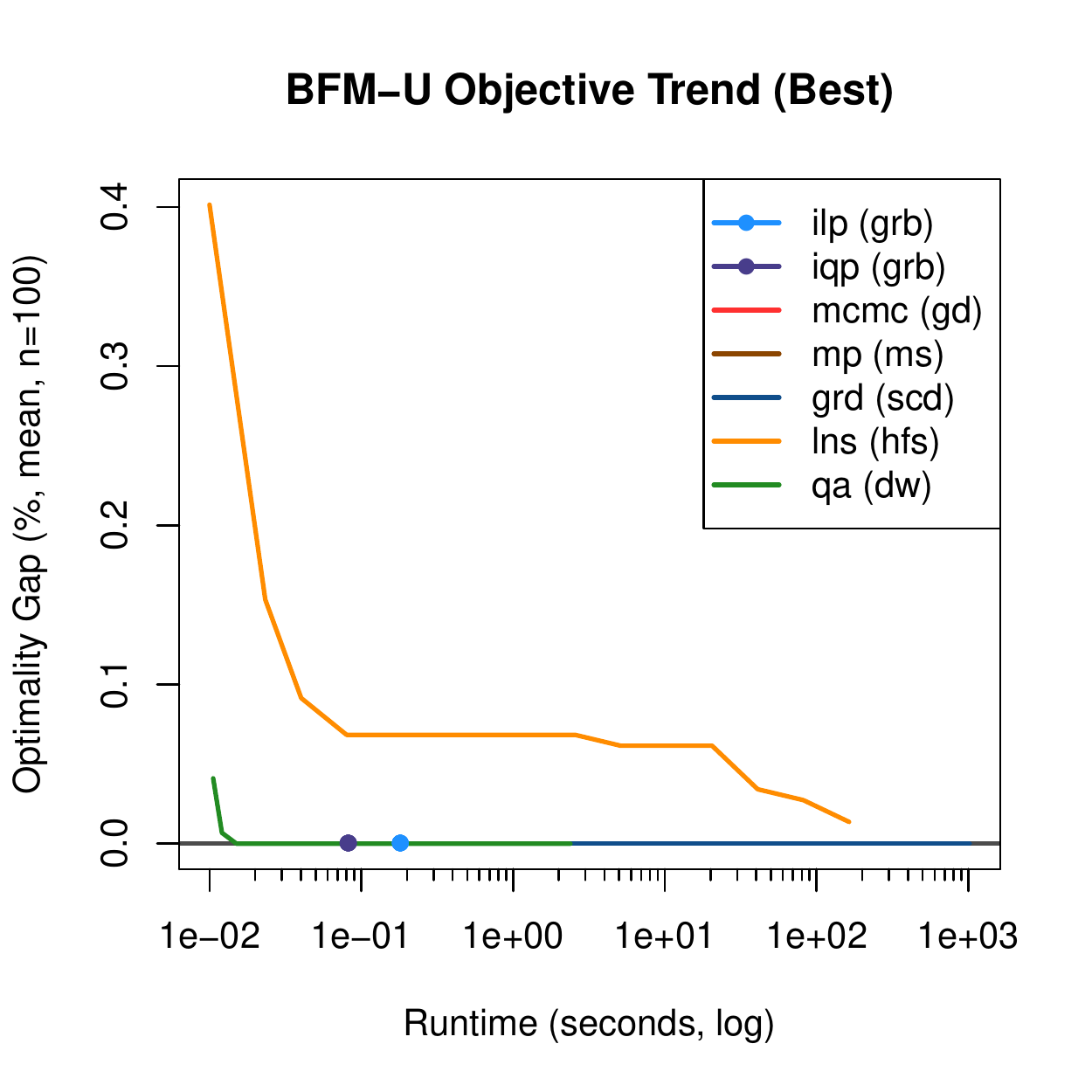}\\
    \vspace{-0.3cm}
    \includegraphics[width=0.43\textwidth]{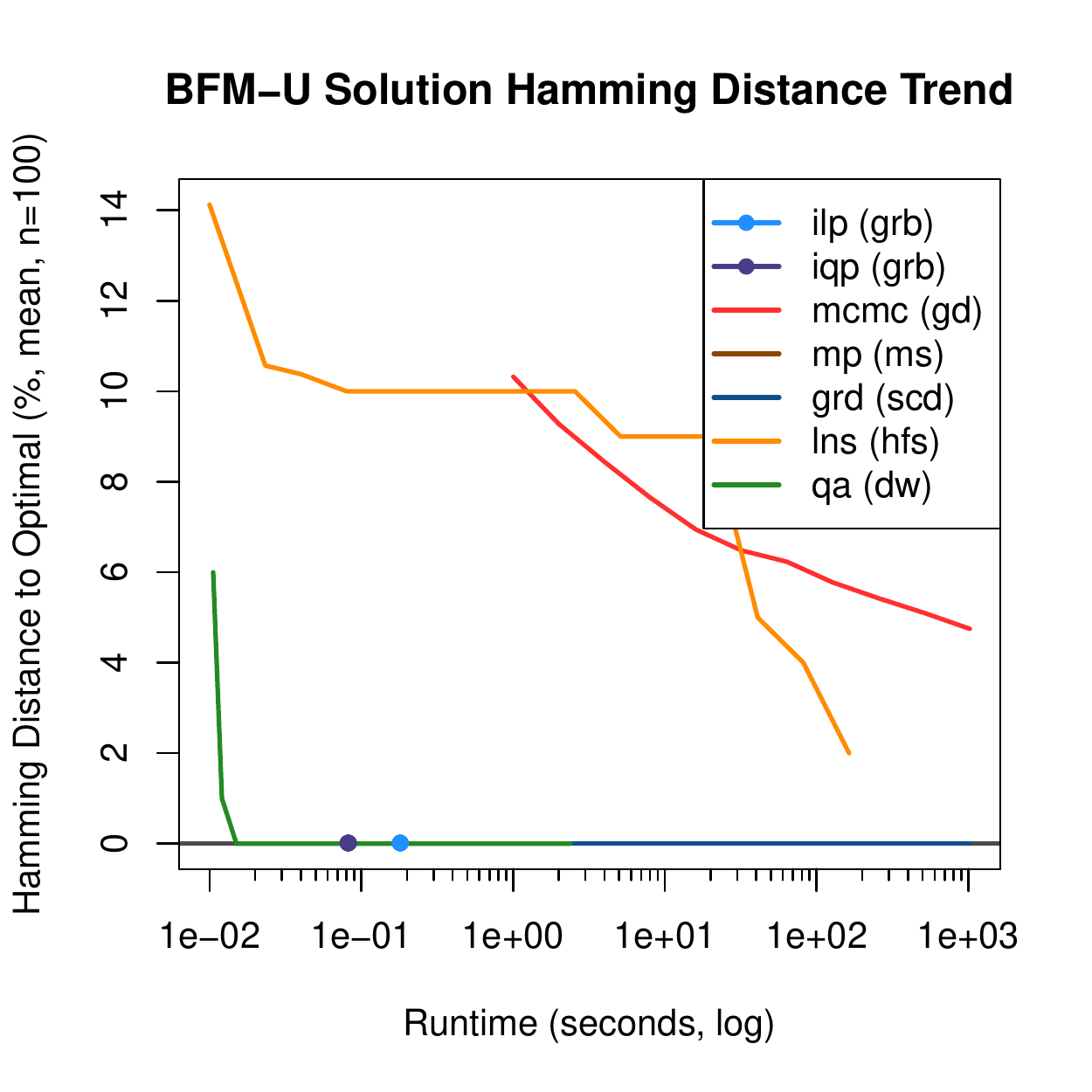}
    \includegraphics[width=0.43\textwidth]{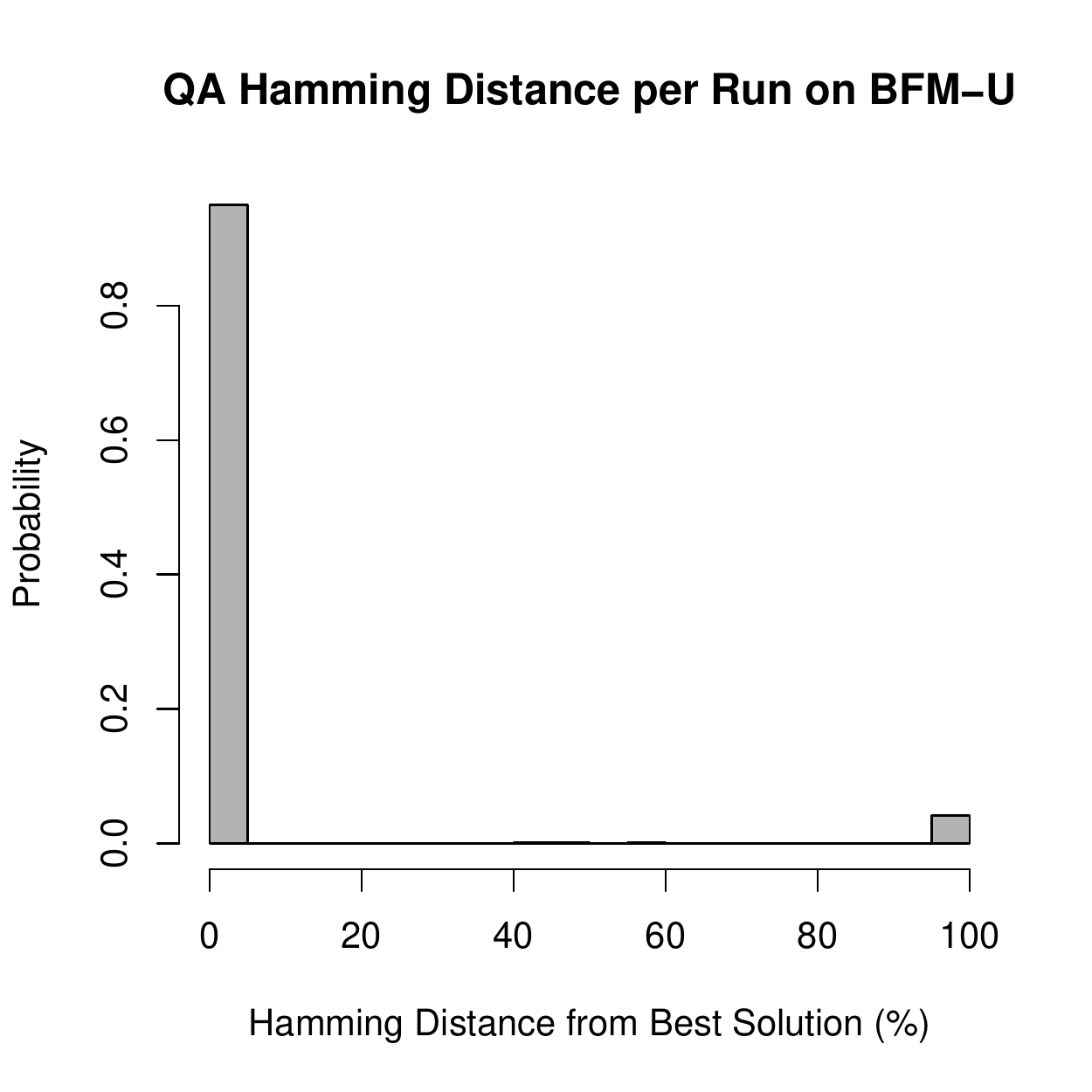}
    \end{center}
    \vspace{-0.5cm}
    \caption{Performance profile (top) and Hamming Distance (bottom) analysis for the Biased Ferromagnet with Uniform Fields instance}
    \label{fig:bfm-u}
\end{figure}

\section{Uniform Fields}
\label{apx:bias}

This appendix presents the results of the uniform-field variants of the BFM, FBFM, and CBFM instances and illustrates how uniform fields improve the performance of all solution methods considered.  Specifically the uniform-field variants replace the bias term, $P(\bm h_i = -1.00) = 0.010$, with the uniform variant $P(\bm h_i = -0.01) = 1.000$.  Throughout this study the field's probability distribution is modified such that there are no zero-value fields (i.e., $P(\bm h_i = 0.00) = 0.000$) and, for consistency with the BFM, FBFM, and CBFM cases presented in Section \ref{sec:structure}, the mean of the fields is selected to be -0.01 (i.e., $\mu_{\bm h} = -0.01$) in all problems considered.

\subsection{The Biased Ferromagnet with Uniform Fields}
\vspace{-0.2cm}
\begin{align}
    & \bm J_{ij} = -1.00 ~\forall i,j \in {\cal E}; \bm h_i = -0.01 ~\forall i \in {\cal N} \label{eq:bfm-u} \tag{BFM-U}
\end{align}

The Biased Ferromagnet with Uniform Fields \eqref{eq:bfm-u} is similar to the BFM case, but all of the linear terms are set identically to $\bm h_i = -0.01$.  All of the solution methods considered here perform well on this BFM-U case (see Figure \ref{fig:bfm-u}).  However, the BFM-U case does appear to reduce both the optimality gap and hamming distance metrics by a factor of two compared to the BFM case.  This suggests that BFM-U is easier than BFM based on the metrics considered by this work.

\begin{figure}[t]
    \begin{center}
    \includegraphics[width=0.43\textwidth]{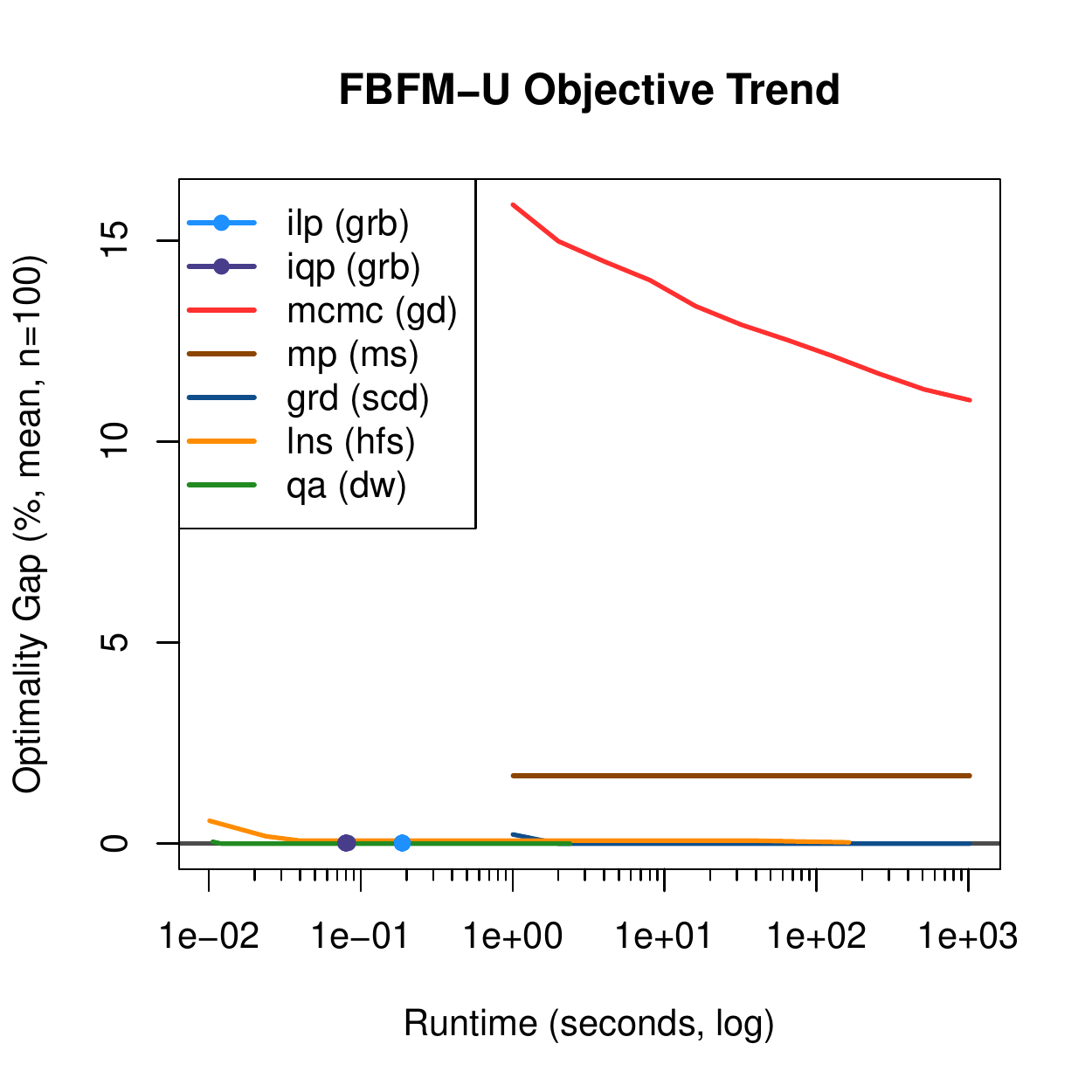}
    \includegraphics[width=0.43\textwidth]{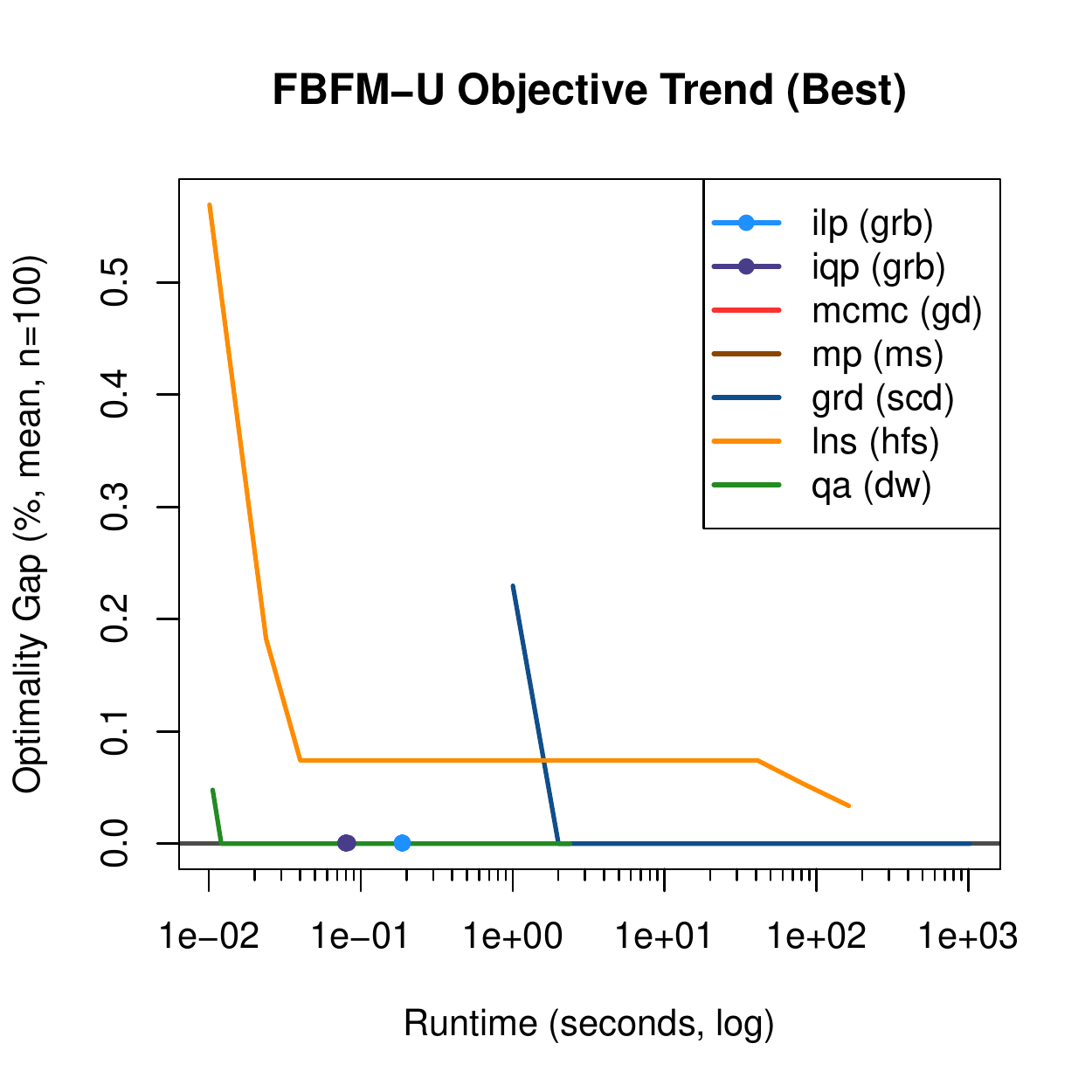}\\
    \vspace{-0.3cm}
    \includegraphics[width=0.43\textwidth]{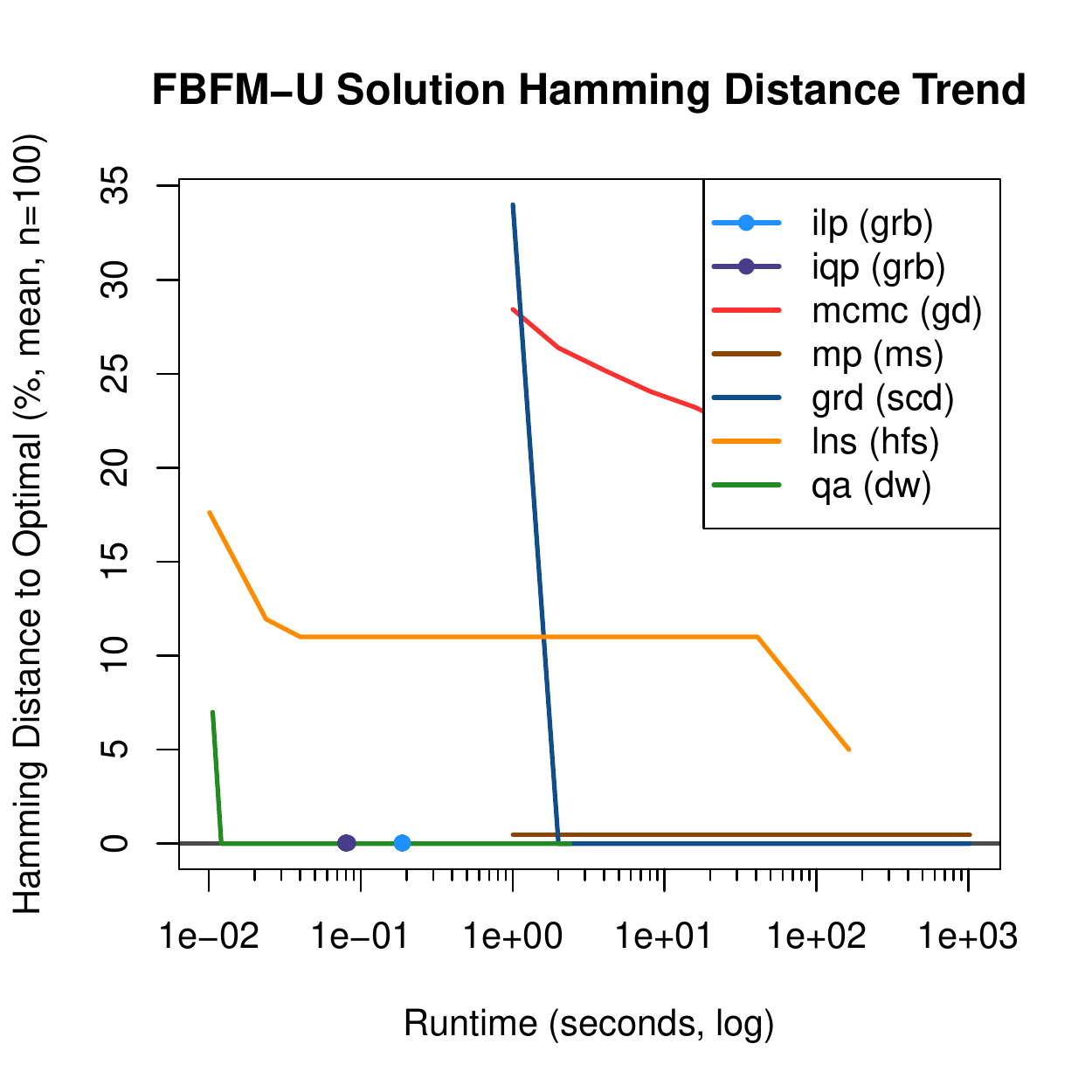}
    \includegraphics[width=0.43\textwidth]{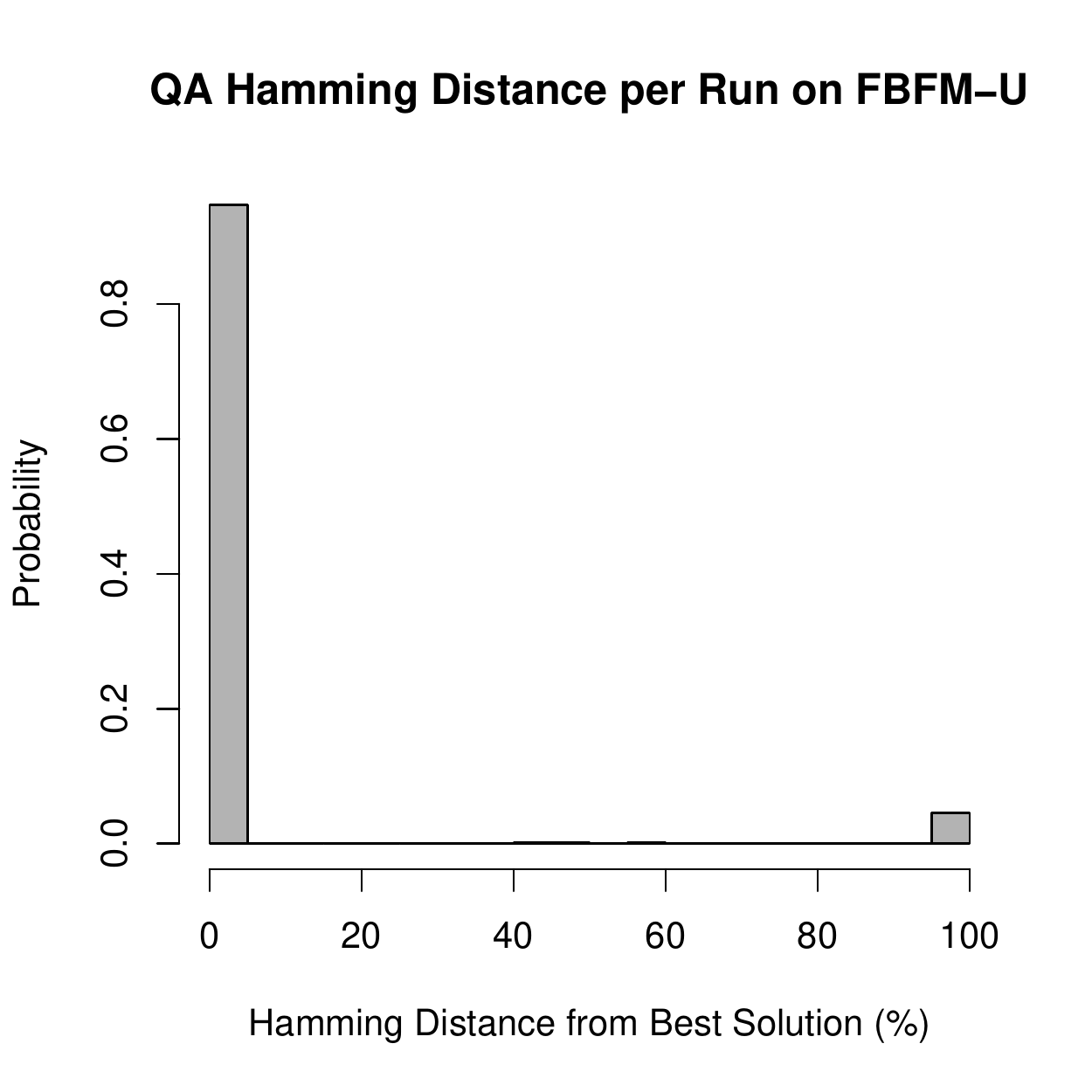}
    \end{center}
    \vspace{-0.5cm}
    \caption{Performance profile (top) and Hamming Distance (bottom) analysis for the Frustrated Biased Ferromagnet with Uniform Fields instance}
    \label{fig:fbfm-u}
\end{figure}

\subsection{The Frustrated Biased Ferromagnet with Uniform Fields}
\vspace{-0.2cm}
\begin{align}
    \label{eq:fbfm-u} \tag{FBFM-U}
    \bm J_{ij} &= -1.00 ~\forall i,j \in {\cal E} \\
    P(\bm h_i = -0.03) = 0.666 &, P(\bm h_i = 0.03) = 0.334 ~\forall i \in {\cal N} \nonumber
\end{align}

The Frustrated Biased Ferromagnet with Uniform Fields \eqref{eq:fbfm-u} is similar to the FBFM case, but two-thirds of the linear terms are set to $\bm h_i = -0.03$ and one-third is set to $\bm h_i = 0.03$.  Although the performance of most of the algorithms on FBFM-U is similar to FBFM (see Figure \ref{fig:fbfm-u}), there are two notable deviations.  The performance of MS and SCD algorithms improves significantly in the FBFM-U case.  This also suggests that the FBFM-U is easier than FBFM based on the metrics considered by this work.

\begin{figure}[t]
    \begin{center}
    \includegraphics[width=0.43\textwidth]{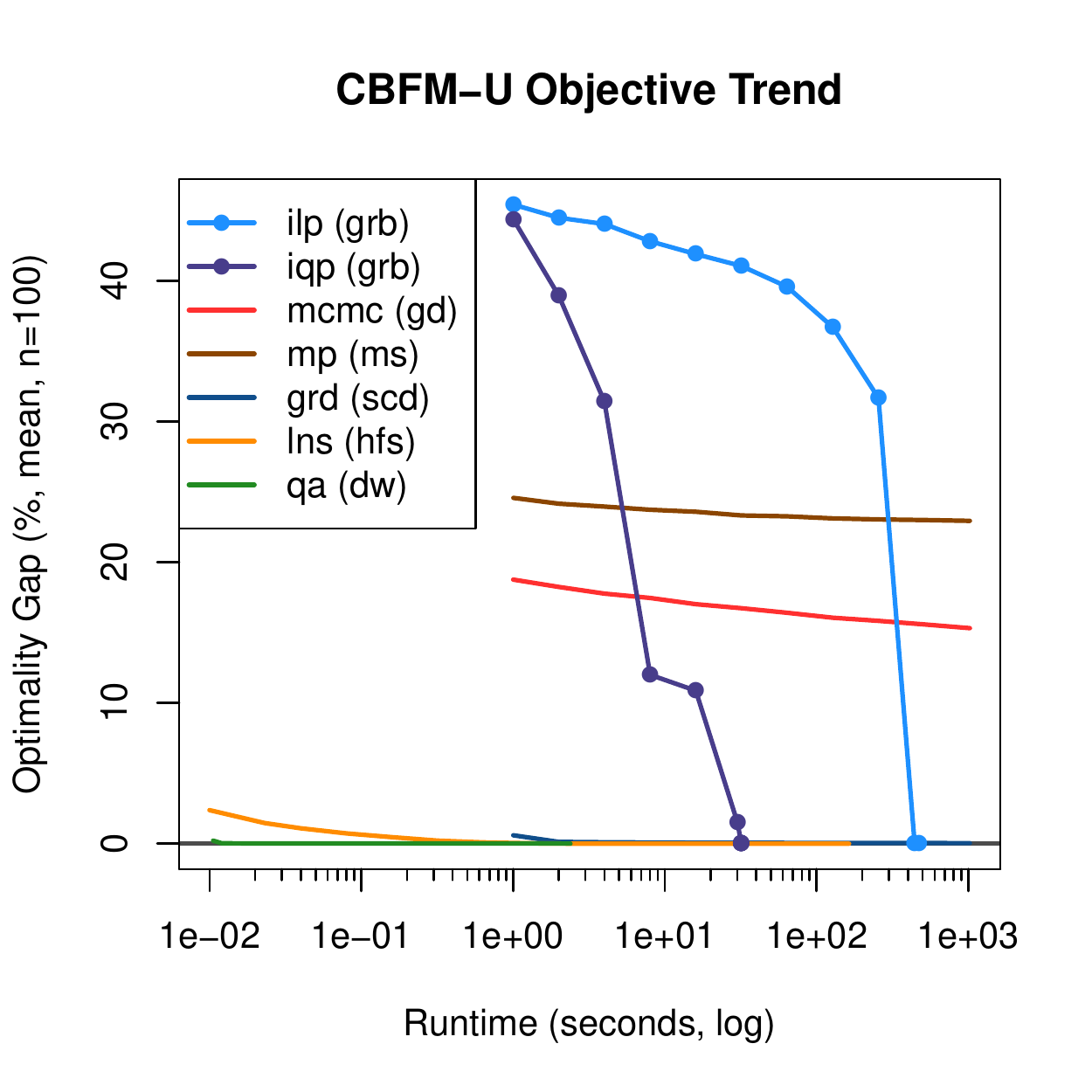}
    \includegraphics[width=0.43\textwidth]{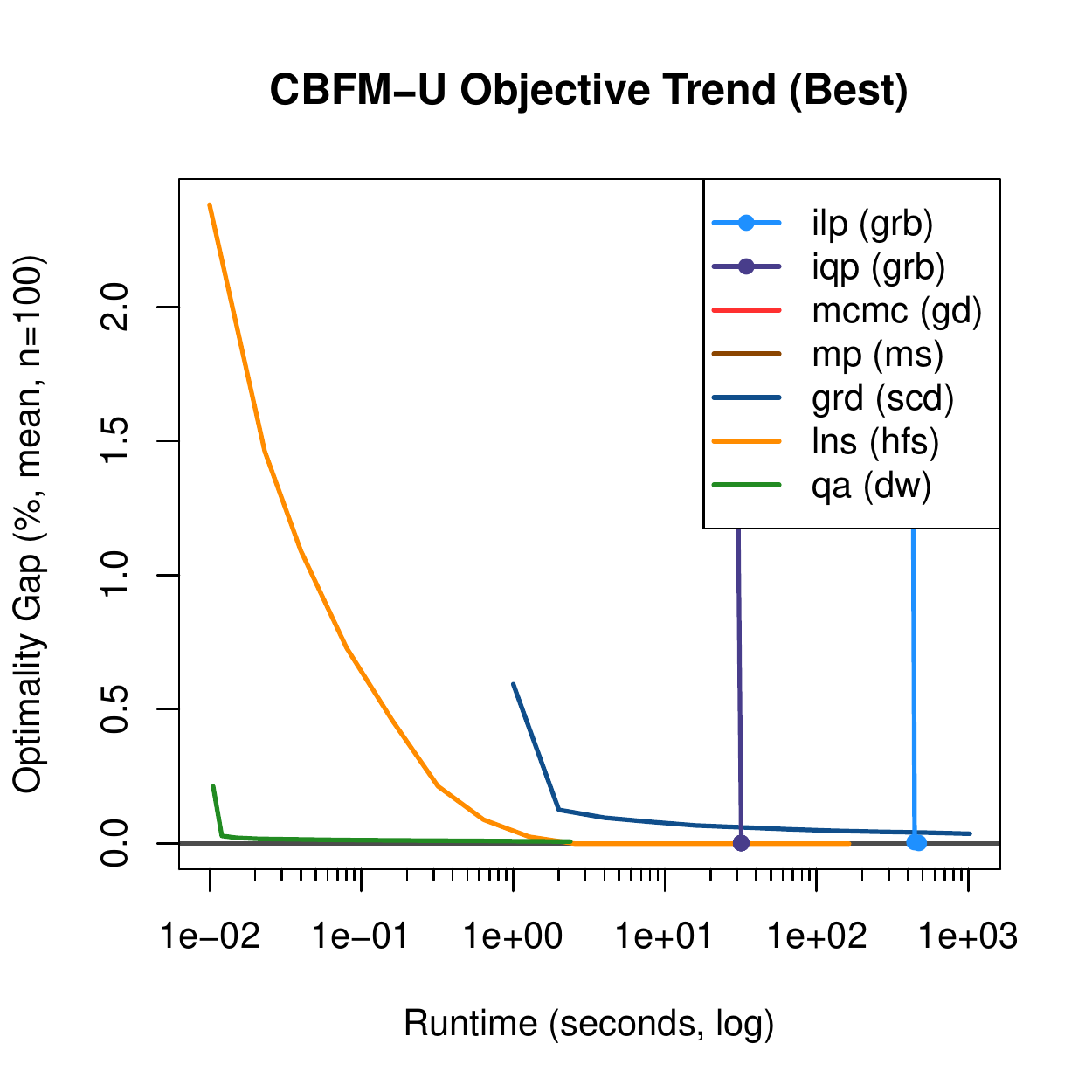}\\
    \vspace{-0.3cm}
    \includegraphics[width=0.43\textwidth]{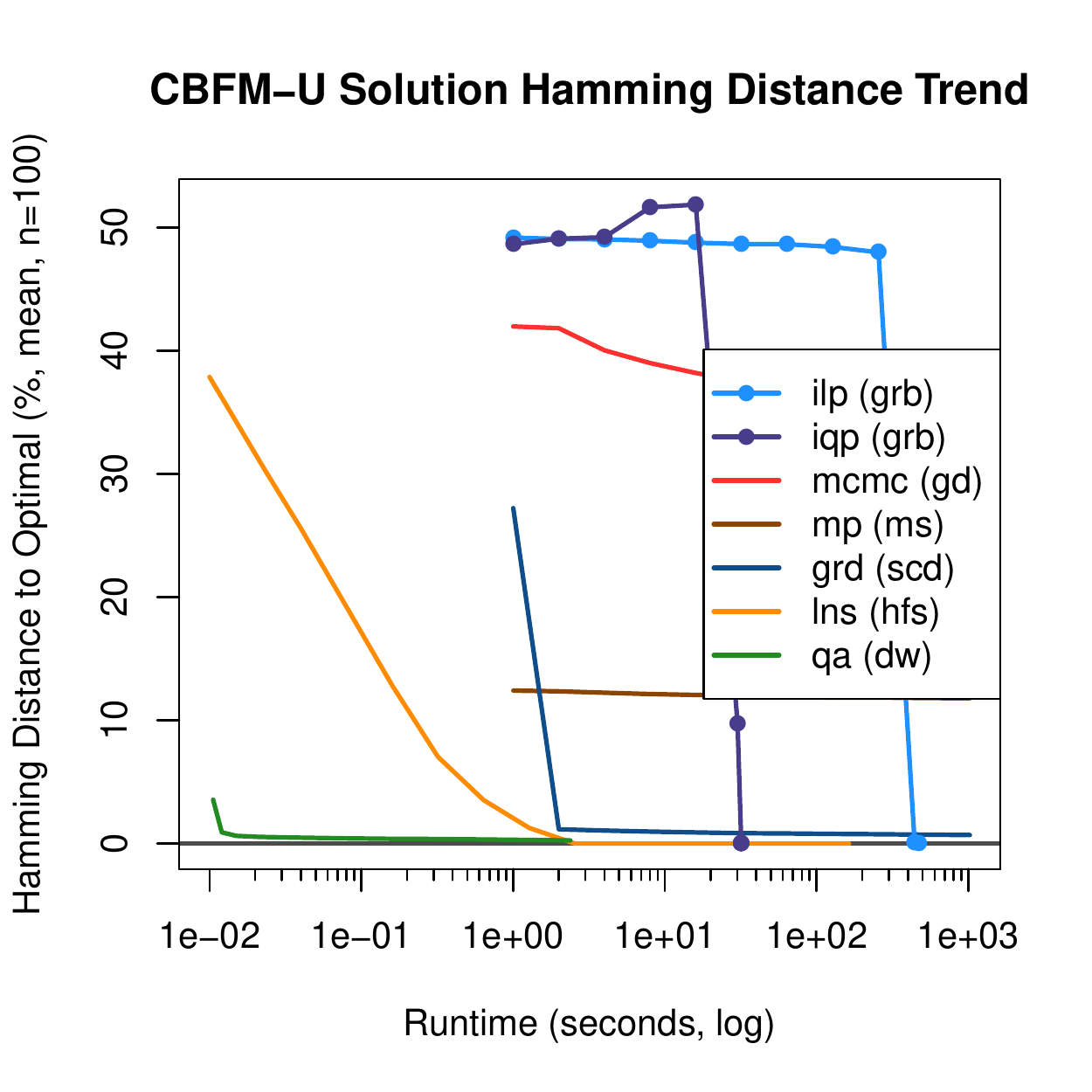}
    \includegraphics[width=0.43\textwidth]{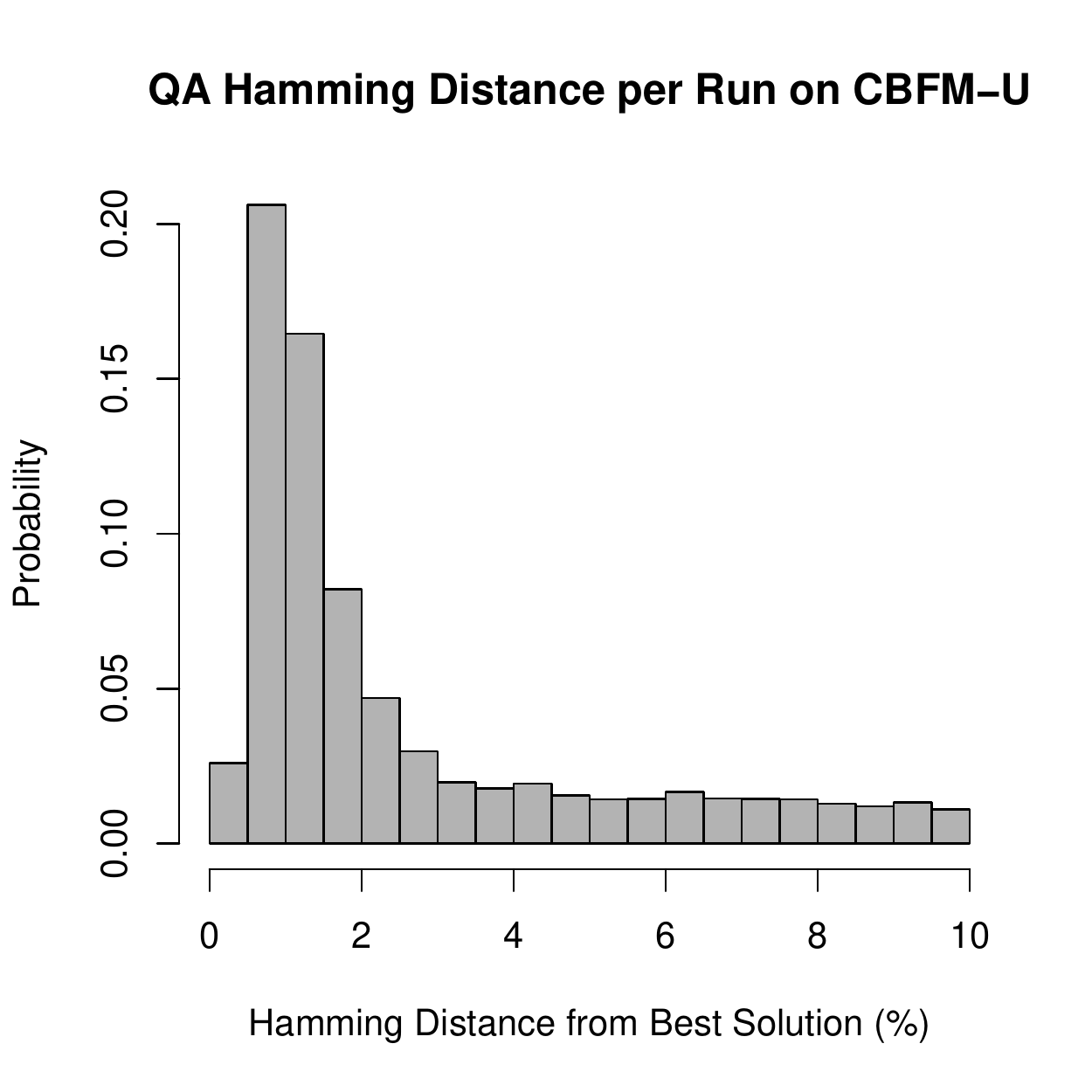}
    \end{center}
    \vspace{-0.5cm}
    \caption{Performance profile (top) and Hamming Distance (bottom) analysis for the Corrupted Biased Ferromagnet with Uniform Fields instance}
    \label{fig:cbfm-u}
\end{figure}

\subsection{The Corrupted Biased Ferromagnet with Uniform Fields}
\vspace{-0.2cm}
\begin{align}
    \label{eq:cbfm-u} \tag{CBFM-U}
    P(\bm J_{ij} = -1.00) = 0.625 &, P(\bm J_{ij} = 0.20) = 0.375 ~\forall i,j \in {\cal E} \\
    P(\bm h_i = -0.03) = 0.666 &, P(\bm h_i = 0.03) = 0.334 ~\forall i \in {\cal N} \nonumber
\end{align}

The Corrupted Biased Ferromagnet with Uniform Fields \eqref{eq:cbfm-u} is similar to the CBFM case, but two-thirds of the linear terms are set to $\bm h_i = -0.03$ and one-third is set to $\bm h_i = 0.03$.  This case exhibits the most variation from the CBFM alternative (see Figure \ref{fig:cbfm-u}).  
The key observations are as follows:
\begin{itemize}
    \item In CBFM-U, QA has a higher probability of finding a near-optimal solution (i.e.,  $>$ 0.50) than CBFM (i.e., $<$ 0.20).  However, it has a lower probability of finding the true-optimal solution (Figure \ref{fig:cbfm-u}, bottom-right).  Due to this effect, QA finds a near-optimal solution to CBFM-U faster than CBFM but never manages to converge to the optimal solution, as it does in CBFM.
    \item The performance of the SCD algorithm improves significantly in the CBFM-U case.  The SCD algorithm is among the best solutions for CBFM-U ($<$ 0.5\% optimality gap), while it has more than a 2\% optimality gap in the CBFM case.
\end{itemize}
Overall, these results suggest that CBFM-U is easier than CBFM based on the metrics considered by this work.  However, the subtle differences in the performance of QA between CBFM and CBFM-U suggest that varying the distribution of the linear terms in the CBFM family of problems could be a useful tool for developing a deeper understanding of how QA responds to different classes of optimization tasks.

\section{Reference Implementations}

\subsection{D-Wave Instance Generator (DWIG)}
The problems considered in this work were generated with the open-source D-Wave Instance Generator tool, which is available at \url{https://github.com/lanl-ansi/dwig}.  DWIG is a command line tool that uses D-Wave's hardware API to identify the topology of a specific D-Wave device and uses that graph for randomized problem generation.  The following list provides the mapping of problems in this paper to the DWIG command line interface:
\begin{verbatim}
CBFM:
  dwig.py cbfm -rgt

CBFM-U:
  dwig.py cbfm -rgt
    -j1-val -1.00 -j1-pr 0.625 -j2-val 0.02 -j2-pr 0.375
    -h1-val -0.03 -h1-pr 0.666 -h2-val 0.03 -h2-pr 0.334

FBFM:
  dwig.py cbfm -rgt
    -j1-val -1.00 -j1-pr 1.000 -j2-val 0.00 -j2-pr 0.000
    -h1-val -1.00 -h1-pr 0.020 -h2-val 1.00 -h2-pr 0.010

FBFM-U:
  dwig.py cbfm -rgt
    -j1-val -1.00 -j1-pr 1.000 -j2-val 0.00 -j2-pr 0.000
    -h1-val -0.03 -h1-pr 0.666 -h2-val 0.03 -h2-pr 0.334

BFM:
  dwig.py cbfm -rgt
    -j1-val -1.00 -j1-pr 1.000 -j2-val 0.00 -j2-pr 0.000
    -h1-val -1.00 -h1-pr 0.010 -h2-val 0.00 -h2-pr 0.000

BFM-U:
  dwig.py cbfm -rgt
    -j1-val -1.00 -j1-pr 1.000 -j2-val 0.00 -j2-pr 0.000
    -h1-val -0.01 -h1-pr 1.000 -h2-val 0.00 -h2-pr 0.000
\end{verbatim}

\clearpage
\subsection{Ising Model Optimization Methods}

The problems considered in this work were solved with the open-source {\em Ising-Solvers scripts} that are available at \url{https://github.com/lanl-ansi/ising-solvers}.  These scripts include a combination of calls to executables, system libraries, and handmade heuristics.  Each script conforms to a standard API for measuring runtime and reporting results.  The following commands were used for each of the solution approaches presented in this work:
\begin{verbatim}
ILP (GRB):
  ilp_gurobi.py -ss -rtl <time_limit> -f <case file>

IQP (GRB):
  iqp_gurobi.py -ss -rtl <time_limit> -f <case file>

MCMC (GD):
  mcmc_gd.py -ss -rtl <time_limit> -f <case file>

MP (MS):
  mp_ms.py -ss -rtl <time_limit> -f <case file>

GRD (SCD):
  grd_scd.jl -s -t <time_limit> -f <case file>

LNS (HFS):
  lns_hfs.py -ss -rtl <time_limit> -f <case file>

QA (DW):
  qa_dwave.py -ss -nr <number of reads> -at 5 -srtr 100 -f <case file>

\end{verbatim}

\clearpage
\bibliographystyle{spmpsci}
\bibliography{references}

\begin{thebibliography}{10}
\providecommand{\url}[1]{{#1}}
\providecommand{\urlprefix}{URL }
\expandafter\ifx\csname urlstyle\endcsname\relax
  \providecommand{\doi}[1]{DOI~\discretionary{}{}{}#1}\else
  \providecommand{\doi}{DOI~\discretionary{}{}{}\begingroup
  \urlstyle{rm}\Url}\fi

\bibitem{aaronson_blog2}
Aaronson, S.: Insert d-wave post here.
\newblock Published online at \url{http://www.scottaaronson.com/blog/?p=3192}
  (2017).
\newblock Accessed: 04/28/2017

\bibitem{Adame_2020}
Adame, J.I., McMahon, P.L.: Inhomogeneous driving in quantum annealers can
  result in orders-of-magnitude improvements in performance.
\newblock Quantum Science and Technology \textbf{5}(3), 035011 (2020).
\newblock \doi{10.1088/2058-9565/ab935a}.
\newblock \urlprefix\url{https://doi.org/10.1088%2F2058-9565%2Fab935a}

\bibitem{albash2018adiabatic}
Albash, T., Lidar, D.A.: Adiabatic quantum computation.
\newblock Reviews of Modern Physics \textbf{90}(1), 015002 (2018)

\bibitem{PhysRevX.8.031016}
Albash, T., Lidar, D.A.: Demonstration of a scaling advantage for a quantum
  annealer over simulated annealing.
\newblock Phys. Rev. X \textbf{8}, 031016 (2018).
\newblock \doi{10.1103/PhysRevX.8.031016}.
\newblock \urlprefix\url{https://link.aps.org/doi/10.1103/PhysRevX.8.031016}

\bibitem{Arute2019}
Arute, F., Arya, K., Babbush, R., Bacon, D., Bardin, J.C., Barends, R., Biswas,
  R., Boixo, S., Brandao, F.G.S.L., et~al.: Quantum supremacy using a
  programmable superconducting processor.
\newblock Nature \textbf{574}(7779), 505--510 (2019).
\newblock \doi{10.1038/s41586-019-1666-5}.
\newblock \urlprefix\url{https://doi.org/10.1038/s41586-019-1666-5}

\bibitem{baccari2018verification}
Baccari, F., Gogolin, C., Wittek, P., Ac{\'\i}n, A.: Verification of quantum
  optimizers.
\newblock arXiv preprint arXiv:1808.01275  (2018)

\bibitem{barahona1982computational}
Barahona, F.: On the computational complexity of ising spin glass models.
\newblock Journal of Physics A: Mathematical and General \textbf{15}(10), 3241
  (1982)

\bibitem{10.3389/fphy.2014.00056}
Bian, Z., Chudak, F., Israel, R., Lackey, B., Macready, W.G., Roy, A.: Discrete
  optimization using quantum annealing on sparse ising models.
\newblock Frontiers in Physics \textbf{2}, 56 (2014).
\newblock \doi{10.3389/fphy.2014.00056}.
\newblock
  \urlprefix\url{http://journal.frontiersin.org/article/10.3389/fphy.2014.00056}

\bibitem{10.3389/fict.2016.00014}
Bian, Z., Chudak, F., Israel, R.B., Lackey, B., Macready, W.G., Roy, A.:
  Mapping constrained optimization problems to quantum annealing with
  application to fault diagnosis.
\newblock Frontiers in ICT \textbf{3}, 14 (2016).
\newblock \doi{10.3389/fict.2016.00014}.
\newblock
  \urlprefix\url{http://journal.frontiersin.org/article/10.3389/fict.2016.00014}

\bibitem{Billionnet2007}
Billionnet, A., Elloumi, S.: Using a mixed integer quadratic programming solver
  for the unconstrained quadratic 0-1 problem.
\newblock Mathematical Programming \textbf{109}(1), 55--68 (2007).
\newblock \doi{10.1007/s10107-005-0637-9}.
\newblock \urlprefix\url{http://dx.doi.org/10.1007/s10107-005-0637-9}

\bibitem{Boixo2014}
Boixo, S., Ronnow, T.F., Isakov, S.V., Wang, Z., Wecker, D., Lidar, D.A.,
  Martinis, J.M., Troyer, M.: Evidence for quantum annealing with more than one
  hundred qubits.
\newblock Nat Phys \textbf{10}(3), 218--224 (2014).
\newblock \urlprefix\url{http://dx.doi.org/10.1038/nphys2900}.
\newblock Article

\bibitem{BOROS2002155}
Boros, E., Hammer, P.L.: Pseudo-boolean optimization.
\newblock Discrete Applied Mathematics \textbf{123}(1), 155 -- 225 (2002).
\newblock \doi{https://doi.org/10.1016/S0166-218X(01)00341-9}.
\newblock
  \urlprefix\url{http://www.sciencedirect.com/science/article/pii/S0166218X01003419}

\bibitem{RevModPhys.39.883}
Brush, S.G.: History of the lenz-ising model.
\newblock Rev. Modern Phys. \textbf{39}, 883--893 (1967).
\newblock \doi{10.1103/RevModPhys.39.883}.
\newblock \urlprefix\url{https://link.aps.org/doi/10.1103/RevModPhys.39.883}

\bibitem{chmielewski2018cloud}
Chmielewski, M., Amini, J., Hudek, K., Kim, J., Mizrahi, J., Monroe, C.,
  Wright, K., Moehring, D.: Cloud-based trapped-ion quantum computing.
\newblock In: APS Meeting Abstracts (2018)

\bibitem{coffrin2016challenges}
Coffrin, C., Nagarajan, H., Bent, R.: {Challenges and Successes of Solving
  Binary Quadratic Programming Benchmarks on the DW2X QPU}.
\newblock Tech. rep., Los Alamos National Laboratory (LANL) (2016)

\bibitem{10.1007/978-3-030-19212-9_11}
Coffrin, C., Nagarajan, H., Bent, R.: Evaluating ising processing units with
  integer programming.
\newblock In: L.M. Rousseau, K.~Stergiou (eds.) Integration of Constraint
  Programming, Artificial Intelligence, and Operations Research, pp. 163--181.
  Springer International Publishing, Cham (2019)

\bibitem{ising_solvers}
Coffrin, C., Pang, Y.: ising-solvers.
\newblock \url{https://github.com/lanl-ansi/ising-solvers} (2019)

\bibitem{coles2018quantum}
Coles, P.J., Eidenbenz, S., Pakin, S., Adedoyin, A., Ambrosiano, J., Anisimov,
  P., Casper, W., Chennupati, G., Coffrin, C., Djidjev, H., et~al.: Quantum
  algorithm implementations for beginners.
\newblock arXiv preprint arXiv:1804.03719  (2018)

\bibitem{ising_frustration}
Cugliandolo, L.F.: Advanced statistical physics: Frustration.
\newblock
  \url{https://www.lpthe.jussieu.fr/~leticia/TEACHING/Master2018/frustration18.pdf}
  (2018)

\bibitem{1306.1202}
Dash, S.: A note on qubo instances defined on chimera graphs.
\newblock arXiv preprint arXiv:1306.1202  (2013).
\newblock \urlprefix\url{https://arxiv.org/abs/1306.1202}

\bibitem{d1985random}
d'Auriac, J.A., Preissmann, M., Rammal, R.: The random field ising model:
  algorithmic complexity and phase transition.
\newblock Journal de Physique Lettres \textbf{46}(5), 173--180 (1985)

\bibitem{PhysRevX.6.031015}
Denchev, V.S., Boixo, S., Isakov, S.V., Ding, N., Babbush, R., Smelyanskiy, V.,
  Martinis, J., Neven, H.: What is the computational value of finite-range
  tunneling?
\newblock Phys. Rev. X \textbf{6}, 031015 (2016).
\newblock \doi{10.1103/PhysRevX.6.031015}.
\newblock \urlprefix\url{https://link.aps.org/doi/10.1103/PhysRevX.6.031015}

\bibitem{dhar1997zero}
Dhar, D., Shukla, P., Sethna, J.P.: Zero-temperature hysteresis in the
  random-field ising model on a bethe lattice.
\newblock Journal of Physics A: Mathematical and General \textbf{30}(15), 5259
  (1997)

\bibitem{ding2015proof}
Ding, J., Sly, A., Sun, N.: Proof of the satisfiability conjecture for large k.
\newblock In: Proceedings of the forty-seventh annual ACM symposium on Theory
  of computing, pp. 59--68. ACM (2015)

\bibitem{eagle2009inferring}
Eagle, N., Pentland, A.S., Lazer, D.: Inferring friendship network structure by
  using mobile phone data.
\newblock Proceedings of the national academy of sciences \textbf{106}(36),
  15274--15278 (2009)

\bibitem{1808.09999}
Fabio L.~Traversa, M.D.V.: Memcomputing integer linear programming (2018).
\newblock \urlprefix\url{https://arxiv.org/abs/1808.09999}

\bibitem{Farhi472}
Farhi, E., Goldstone, J., Gutmann, S., Lapan, J., Lundgren, A., Preda, D.: A
  quantum adiabatic evolution algorithm applied to random instances of an
  np-complete problem.
\newblock Science \textbf{292}(5516), 472--475 (2001).
\newblock \doi{10.1126/science.1057726}.
\newblock \urlprefix\url{http://science.sciencemag.org/content/292/5516/472}

\bibitem{quant-ph-0001106}
Farhi, E., Goldstone, J., Gutmann, S., Sipser, M.: Quantum computation by
  adiabatic evolution (2018).
\newblock \urlprefix\url{https://arxiv.org/abs/quant-ph/0001106}

\bibitem{Feynman1982-FEYSPW}
Feynman, R.P.: Simulating physics with computers.
\newblock International Journal of Theoretical Physics \textbf{21}(6), 467--488
  (1982)

\bibitem{fossorier1999reduced}
Fossorier, M.P., Mihaljevic, M., Imai, H.: Reduced complexity iterative
  decoding of low-density parity check codes based on belief propagation.
\newblock IEEE Transactions on communications \textbf{47}(5), 673--680 (1999)

\bibitem{fuhitsu_da}
Fujitsu: Digital annealer.
\newblock Published online at
  \url{http://www.fujitsu.com/global/digitalannealer/} (2018).
\newblock Accessed: 02/26/2019

\bibitem{gallavotti2013statistical}
Gallavotti, G.: Statistical mechanics: A short treatise.
\newblock Springer Science \& Business Media (2013)

\bibitem{glauber1963time}
Glauber, R.J.: Time-dependent statistics of the ising model.
\newblock Journal of mathematical physics \textbf{4}(2), 294--307 (1963)

\bibitem{grover}
Grover, L.K.: A fast quantum mechanical algorithm for database search.
\newblock In: Proceedings of the twenty-eighth annual ACM symposium on Theory
  of computing, pp. 212--219. ACM (1996)

\bibitem{gurobi}
{Gurobi Optimization, Inc.}: Gurobi optimizer reference manual.
\newblock Published online at \url{http://www.gurobi.com} (2014)

\bibitem{hamerly2019experimental}
Hamerly, R., Inagaki, T., McMahon, P.L., Venturelli, D., Marandi, A., Onodera,
  T., Ng, E., Langrock, C., Inaba, K., Honjo, T., et~al.: Experimental
  investigation of performance differences between coherent ising machines and
  a quantum annealer.
\newblock Science advances \textbf{5}(5), eaau0823 (2019)

\bibitem{Hamze:2004:FT:1036843.1036873}
Hamze, F., de~Freitas, N.: From fields to trees.
\newblock In: Proceedings of the 20th Conference on Uncertainty in Artificial
  Intelligence, UAI '04, pp. 243--250. AUAI Press, Arlington, Virginia, United
  States (2004).
\newblock \urlprefix\url{http://dl.acm.org/citation.cfm?id=1036843.1036873}

\bibitem{Haribara2016}
Haribara, Y., Utsunomiya, S., Yamamoto, Y.: A Coherent Ising Machine for
  MAX-CUT Problems: Performance Evaluation against Semidefinite Programming and
  Simulated Annealing, pp. 251--262.
\newblock Springer Japan, Tokyo (2016).
\newblock \doi{10.1007/978-4-431-55756-2\_12}.
\newblock \urlprefix\url{http://dx.doi.org/10.1007/978-4-431-55756-2\_12}

\bibitem{hopfield1982neural}
Hopfield, J.J.: Neural networks and physical systems with emergent collective
  computational abilities.
\newblock Proceedings of the national academy of sciences \textbf{79}(8),
  2554--2558 (1982)

\bibitem{Inagaki603}
Inagaki, T., Haribara, Y., Igarashi, K., Sonobe, T., Tamate, S., Honjo, T.,
  Marandi, A., McMahon, P.L., Umeki, T., Enbutsu, K., Tadanaga, O., Takenouchi,
  H., Aihara, K., Kawarabayashi, K.i., Inoue, K., Utsunomiya, S., Takesue, H.:
  A coherent ising machine for 2000-node optimization problems.
\newblock Science \textbf{354}(6312), 603--606 (2016).
\newblock \doi{10.1126/science.aah4243}.
\newblock \urlprefix\url{http://science.sciencemag.org/content/354/6312/603}

\bibitem{ibm_quantum}
{International Business Machines Corporation}: Ibm building first universal
  quantum computers for business and science.
\newblock Published online at
  \url{https://www-03.ibm.com/press/us/en/pressrelease/51740.wss} (2017).
\newblock Accessed: 04/28/2017

\bibitem{ISAKOV2015265}
Isakov, S., Zintchenko, I., Rønnow, T., Troyer, M.: Optimised simulated
  annealing for ising spin glasses.
\newblock Computer Physics Communications \textbf{192}, 265 -- 271 (2015).
\newblock \doi{https://doi.org/10.1016/j.cpc.2015.02.015}.
\newblock
  \urlprefix\url{http://www.sciencedirect.com/science/article/pii/S0010465515000727}

\bibitem{johnson2011quantum}
Johnson, M.W., Amin, M.H., Gildert, S., Lanting, T., Hamze, F., Dickson, N.,
  Harris, R., Berkley, A.J., Johansson, J., Bunyk, P., et~al.: Quantum
  annealing with manufactured spins.
\newblock Nature \textbf{473}(7346), 194--198 (2011)

\bibitem{junger2019performance}
J{\"u}nger, M., Lobe, E., Mutzel, P., Reinelt, G., Rendl, F., Rinaldi, G.,
  Stollenwerk, T.: Performance of a quantum annealer for ising ground state
  computations on chimera graphs.
\newblock arXiv preprint arXiv:1904.11965  (2019)

\bibitem{PhysRevE.58.5355}
Kadowaki, T., Nishimori, H.: Quantum annealing in the transverse ising model.
\newblock Phys. Rev. E \textbf{58}, 5355--5363 (1998).
\newblock \doi{10.1103/PhysRevE.58.5355}.
\newblock \urlprefix\url{https://link.aps.org/doi/10.1103/PhysRevE.58.5355}

\bibitem{kalinin2018global}
Kalinin, K.P., Berloff, N.G.: Global optimization of spin hamiltonians with
  gain-dissipative systems.
\newblock Scientific reports \textbf{8}(1), 1--9 (2018)

\bibitem{7738704}
Kielpinski, D., Bose, R., Pelc, J., Vaerenbergh, T.V., Mendoza, G., Tezak, N.,
  Beausoleil, R.G.: Information processing with large-scale optical integrated
  circuits.
\newblock In: 2016 IEEE International Conference on Rebooting Computing (ICRC),
  pp. 1--4 (2016).
\newblock \doi{10.1109/ICRC.2016.7738704}

\bibitem{king2015performance}
King, A.D., Lanting, T., Harris, R.: Performance of a quantum annealer on
  range-limited constraint satisfaction problems.
\newblock arXiv preprint arXiv:1502.02098  (2015)

\bibitem{1701.04579}
King, J., Yarkoni, S., Raymond, J., Ozfidan, I., King, A.D., Nevisi, M.M.,
  Hilton, J.P., McGeoch, C.C.: Quantum annealing amid local ruggedness and
  global frustration (2017).
\newblock \urlprefix\url{https://arxiv.org/abs/1701.04579}

\bibitem{PhysRevA.96.042322}
Lanting, T., King, A.D., Evert, B., Hoskinson, E.: Experimental demonstration
  of perturbative anticrossing mitigation using nonuniform driver hamiltonians.
\newblock Phys. Rev. A \textbf{96}, 042322 (2017).
\newblock \doi{10.1103/PhysRevA.96.042322}.
\newblock \urlprefix\url{https://link.aps.org/doi/10.1103/PhysRevA.96.042322}

\bibitem{leleu2019destabilization}
Leleu, T., Yamamoto, Y., McMahon, P.L., Aihara, K.: Destabilization of local
  minima in analog spin systems by correction of amplitude heterogeneity.
\newblock Physical review letters \textbf{122}(4), 040607 (2019)

\bibitem{lokhov2018optimal}
Lokhov, A.Y., Vuffray, M., Misra, S., Chertkov, M.: Optimal structure and
  parameter learning of ising models.
\newblock Science advances \textbf{4}(3), e1700791 (2018)

\bibitem{10.3389/fphy.2014.00005}
Lucas, A.: Ising formulations of many np problems.
\newblock Frontiers in Physics \textbf{2}, 5 (2014).
\newblock \doi{10.3389/fphy.2014.00005}.
\newblock
  \urlprefix\url{http://journal.frontiersin.org/article/10.3389/fphy.2014.00005}

\bibitem{PhysRevA.94.022337}
Mandr\`a, S., Zhu, Z., Wang, W., Perdomo-Ortiz, A., Katzgraber, H.G.: Strengths
  and weaknesses of weak-strong cluster problems: A detailed overview of
  state-of-the-art classical heuristics versus quantum approaches.
\newblock Phys. Rev. A \textbf{94}, 022337 (2016).
\newblock \doi{10.1103/PhysRevA.94.022337}.
\newblock \urlprefix\url{https://link.aps.org/doi/10.1103/PhysRevA.94.022337}

\bibitem{marbach2012wisdom}
Marbach, D., Costello, J.C., K{\"u}ffner, R., Vega, N.M., Prill, R.J., Camacho,
  D.M., Allison, K.R., Aderhold, A., Bonneau, R., Chen, Y., et~al.: Wisdom of
  crowds for robust gene network inference.
\newblock Nature methods \textbf{9}(8), 796 (2012)

\bibitem{PhysRevApplied.11.044083}
Marshall, J., Venturelli, D., Hen, I., Rieffel, E.G.: Power of pausing:
  Advancing understanding of thermalization in experimental quantum annealers.
\newblock Phys. Rev. Applied \textbf{11}, 044083 (2019).
\newblock \doi{10.1103/PhysRevApplied.11.044083}.
\newblock
  \urlprefix\url{https://link.aps.org/doi/10.1103/PhysRevApplied.11.044083}

\bibitem{dwave_ocp}
McGeoch, C.C., King, J., Nevisi, M.M., Yarkoni, S., Hilton, J.: Optimization
  with clause problems.
\newblock Published online at
  \url{https://www.dwavesys.com/sites/default/files/14-1001A_tr_Optimization_with_Clause_Problems.pdf}
  (2017).
\newblock Accessed: 02/10/2020

\bibitem{McGeoch:2013:EEA:2482767.2482797}
McGeoch, C.C., Wang, C.: Experimental evaluation of an adiabiatic quantum
  system for combinatorial optimization.
\newblock In: Proceedings of the ACM International Conference on Computing
  Frontiers, CF '13, pp. 23:1--23:11. ACM, New York, NY, USA (2013).
\newblock \doi{10.1145/2482767.2482797}.
\newblock \urlprefix\url{http://doi.acm.org/10.1145/2482767.2482797}

\bibitem{mcmahon2016fully}
McMahon, P.L., Marandi, A., Haribara, Y., Hamerly, R., Langrock, C., Tamate,
  S., Inagaki, T., Takesue, H., Utsunomiya, S., Aihara, K., et~al.: A
  fully-programmable 100-spin coherent ising machine with all-to-all
  connections.
\newblock Science p. aah5178 (2016)

\bibitem{mezard2009information}
Mezard, M., Mezard, M., Montanari, A.: Information, physics, and computation.
\newblock Oxford University Press (2009)

\bibitem{mezard1985microstructure}
M{\'e}zard, M., Virasoro, M.A.: The microstructure of ultrametricity.
\newblock Journal de Physique \textbf{46}(8), 1293--1307 (1985)

\bibitem{45919}
Mohseni, M., Read, P., Neven, H., Boixo, S., Denchev, V., Babbush, R., Fowler,
  A., Smelyanskiy, V., Martinis, J.: Commercialize quantum technologies in five
  years.
\newblock Nature \textbf{543}, 171–174 (2017).
\newblock
  \urlprefix\url{http://www.nature.com/news/commercialize-quantum-technologies-in-five-years-1.21583}

\bibitem{morcos2011direct}
Morcos, F., Pagnani, A., Lunt, B., Bertolino, A., Marks, D.S., Sander, C.,
  Zecchina, R., Onuchic, J.N., Hwa, T., Weigt, M.: Direct-coupling analysis of
  residue coevolution captures native contacts across many protein families.
\newblock Proceedings of the National Academy of Sciences \textbf{108}(49),
  E1293--E1301 (2011)

\bibitem{panjwani1995markov}
Panjwani, D.K., Healey, G.: Markov random field models for unsupervised
  segmentation of textured color images.
\newblock IEEE Transactions on pattern analysis and machine intelligence
  \textbf{17}(10), 939--954 (1995)

\bibitem{1604.00319}
Parekh, O., Wendt, J., Shulenburger, L., Landahl, A., Moussa, J., Aidun, J.:
  Benchmarking adiabatic quantum optimization for complex network analysis
  (2015).
\newblock \urlprefix\url{https://arxiv.org/abs/1604.00319}

\bibitem{ibm_blog}
Puget, J.F.: D-wave vs cplex comparison. part 2: Qubo.
\newblock Published online (2013).
\newblock Accessed: 11/28/2018

\bibitem{Rieffel2015}
Rieffel, E.G., Venturelli, D., O'Gorman, B., Do, M.B., Prystay, E.M.,
  Smelyanskiy, V.N.: A case study in programming a quantum annealer for hard
  operational planning problems.
\newblock Quantum Information Processing \textbf{14}(1), 1--36 (2015).
\newblock \doi{10.1007/s11128-014-0892-x}.
\newblock \urlprefix\url{http://dx.doi.org/10.1007/s11128-014-0892-x}

\bibitem{schneidman2006weak}
Schneidman, E., Berry~II, M.J., Segev, R., Bialek, W.: Weak pairwise
  correlations imply strongly correlated network states in a neural population.
\newblock Nature \textbf{440}(7087), 1007 (2006)

\bibitem{HFS_impl_2017}
Selby, A.: Qubo-chimera.
\newblock \url{https://github.com/alex1770/QUBO-Chimera} (2013)

\bibitem{1409.3934}
Selby, A.: Efficient subgraph-based sampling of ising-type models with
  frustration (2014).
\newblock \urlprefix\url{https://arxiv.org/abs/1409.3934}

\bibitem{shor1994algorithms}
Shor, P.W.: Algorithms for quantum computation: Discrete logarithms and
  factoring.
\newblock In: Proceedings 35th annual symposium on foundations of computer
  science, pp. 124--134. Ieee (1994)

\bibitem{1506.08479}
Venturelli, D., Marchand, D.J.J., Rojo, G.: Quantum annealing implementation of
  job-shop scheduling (2015).
\newblock \urlprefix\url{https://arxiv.org/abs/1506.08479}

\bibitem{vuffray2014cavity}
Vuffray, M.: The cavity method in coding theory.
\newblock Tech. rep., EPFL (2014)

\bibitem{NIPS2016_6375}
Vuffray, M., Misra, S., Lokhov, A., Chertkov, M.: Interaction screening:
  Efficient and sample-optimal learning of ising models.
\newblock In: D.D. Lee, M.~Sugiyama, U.V. Luxburg, I.~Guyon, R.~Garnett (eds.)
  Advances in Neural Information Processing Systems 29, pp. 2595--2603. Curran
  Associates, Inc. (2016)

\bibitem{vuffray2019efficient}
Vuffray, M., Misra, S., Lokhov, A.Y.: Efficient learning of discrete graphical
  models.
\newblock arXiv preprint arXiv:1902.00600  (2019)

\bibitem{7063111}
Yamaoka, M., Yoshimura, C., Hayashi, M., Okuyama, T., Aoki, H., Mizuno, H.:
  24.3 20k-spin ising chip for combinational optimization problem with cmos
  annealing.
\newblock In: 2015 IEEE International Solid-State Circuits Conference - (ISSCC)
  Digest of Technical Papers, pp. 1--3 (2015).
\newblock \doi{10.1109/ISSCC.2015.7063111}

\bibitem{6662276}
Yoshimura, C., Yamaoka, M., Aoki, H., Mizuno, H.: Spatial computing
  architecture using randomness of memory cell stability under voltage control.
\newblock In: 2013 European Conference on Circuit Theory and Design (ECCTD),
  pp. 1--4 (2013).
\newblock \doi{10.1109/ECCTD.2013.6662276}

\end{thebibliography}
LA-UR-20-22733

\end{document}